\documentclass[11pt]{amsart}
\usepackage{amsmath,amssymb,amsthm,mathtools}
\usepackage[top=3cm,bottom=2cm,left=3cm,right=3cm,marginparwidth=1.75cm]{geometry}
\usepackage[shortlabels]{enumitem}
\usepackage[colorlinks=true, allcolors=blue]{hyperref}
\usepackage[capitalise]{cleveref}
\usepackage{tikz}
\usepackage[useregional]{datetime2}
\usepackage{cancel}
\usepackage{xcolor}
\usepackage{comment}

\newtheorem{theorem}{Theorem}[section]
\newtheorem{proposition}[theorem]{Proposition}
\newtheorem{lemma}[theorem]{Lemma}
\newtheorem{corollary}[theorem]{Corollary}
\theoremstyle{definition}
\newtheorem{definition}[theorem]{Definition}
\newtheorem{remark}[theorem]{Remark}

\usepackage{comment}

\numberwithin{equation}{section}

\newcommand{\R}{\mathbb{R}}

\newcommand{\cG}{\mathcal{G}}
\newcommand{\cV}{\mathcal{V}}

\DeclareMathOperator{\id}{Id}

\DeclareMathOperator{\tr}{tr}

\DeclareMathOperator{\codim}{codim}

\DeclareMathOperator{\vol}{vol}
\DeclareMathOperator{\crit}{crit}
\DeclareMathOperator{\sym}{Sym}
\DeclareMathOperator{\dive}{div}

\usepackage{environ}
\newtheorem{theoremA}{Theorem}

\crefname{theoremA}{Theorem}{Theorems}
\newtheorem{theoremB}{Theorem}

\crefname{theoremB}{Theorem}{Theorems}

\title{Flexibility of Laplace--Beltrami eigenfunctions under metric variation}
\date{}
\author{Agrachev A. Andrei, Da Ponte Nicola, Lerario Antonio, Testa Matteo}
\begin{document}
\begin{abstract}
We study the surjectivity of the differential of the Laplace--Beltrami operator with respect to the metric parameter. We use the obtained results to prove that every small deformation of a nonconstant eigenfunction can be realized with a small perturbation of the metric. As a consequence, we prove that every generic property of smooth functions which is invariant under multiplication by a nonzero scalar is generic for eigenfunctions of the Laplace--Beltrami operator, extending a result by Uhlenbeck.
\end{abstract}
\maketitle

The Laplace--Beltrami operator is a classical differential operator associated with a Riemannian manifold. Its spectral theory connects analysis, geometry and topology, and a natural question is to understand which properties of its eigenfunctions are stable, or generic, when the metric is allowed to vary. In \cite{Uhlenbeck}, Uhlenbeck proved that on a smooth compact connected manifold of dimension at least two, the generic\footnote{Here and throughout the paper, generic means belonging to a residual set, that is, to a set containing a countable intersection of open dense sets.}
 $C^k$ metric has simple eigenvalues and eigenfunctions which are Morse and transverse to zero. The aim of this paper is to prove a general form of this statement in the smooth category: informally, we will show that every generic property of smooth functions, in the sense of the Fr\'echet topology, is also generic for Laplace--Beltrami eigenfunctions. 

We explain now the main ideas and statements in more details. Let $M$ be a compact, connected, smooth manifold without boundary and of dimension at least two. We denote by $\mathcal G_+$ the space of smooth Riemannian metrics on $M$. The main analytic input is a submersivity result for the Laplace--Beltrami operator with respect to the metric variable, that we explain first. Consider the smooth map
\[
\Delta:\mathcal G_+\times C^\infty(M)\longrightarrow C^\infty(M),
\qquad(g,u)\longmapsto \Delta_g u,
\]
and denote its differential in the metric direction by
\[
\dot\Delta_g[h](u)\coloneqq\left.\frac{d}{dt}\right|_{t=0}\Delta_{g+th}u.
\]
 We introduce the open dense set $\mathcal N\subset \mathcal G_+\times C^\infty(M)$ defined as
\[
\mathcal N=\begin{cases}\{(g,u):u\text{ is Morse and }\crit(u)\cap(\Delta_g u)^{-1}(0)=\varnothing\},& \dim M=2,\\[2mm]\{(g,u):u\text{ is Morse}\},& \dim M>2.
\end{cases}
\]
In \cref{section N} we will prove that for every $(g_0,u_0)\in\mathcal N$ there is a neighborhood $V$ of $(g_0,u_0)$ such that the family
\[
\dot\Delta:V\times\mathcal G\longrightarrow C^\infty(M),
\qquad
((g,u),h)\longmapsto \dot\Delta_g[h](u),
\]
is surjective and admits a smooth tame family of right inverses. See \cref{laplacian use}. Here we use the word ``tame'' in the sense of \cite{Hamilton}. Notice that this statement does not require $u$ to be an eigenfunction. In dimension at least three it therefore says, in particular, that for every Morse function $u$ the map $g\mapsto\Delta_gu$ is locally open, see \cref{dimension 3 cor}. The one-dimensional case is genuinely different: there the metric differential of the Laplace--Beltrami operator is never surjective for a nonconstant function. See \cref{no 1}.

We then apply \cref{laplacian use} to the \emph{eigenpair space}. Consider the space of simple eigenpairs
\[
Q'=\left\{(g,\lambda,u)\in\mathcal G_+\times\R_{>0} \times C^\infty(M):(\Delta_g+\lambda)u=0,\ m_\lambda(g)=1,\ u\neq0\right\},
\]
where $m_\lambda=\mathrm{dim}(\ker(\Delta_g+\lambda))$. We further define
\[
Q=\{(g,\lambda,u)\in Q':(g,u)\in\mathcal N\}.
\]
(There is a technical reason to work with this set, that will become more transparent below.)
Using Uhlenbeck's genericity theorem at every sufficiently high regularity, we show that $Q$ is open and dense in $Q'$. The submersivity result above then gives the main ``flexibility'' statement.

\begin{theoremA}\label{theorem A}
Let $M$ be a compact connected smooth manifold of dimension at least two. The projection
\[
\pi:Q\longrightarrow C^\infty(M),\qquad(g,\lambda,u)\longmapsto u,
\]
is an open map.
\end{theoremA}

Equivalently, if $(g_0,\lambda_0,u_0)\in Q$, then every sufficiently small smooth deformation of $u_0$ is again an eigenfunction of a nearby metric, for a nearby eigenvalue. This is the sense in which eigenfunctions are flexible under metric variation, explaining the title of the paper. The statement is local on the eigenpair space, but the subset on which it holds is open and dense.

We explain now the genericity statement generalizing Uhlenbeck's result in the above sense. To proceed, since eigenfunctions are naturally defined up to multiplication by a nonzero scalar, it will be convenient to formulate the result using the projectivized space $\mathbb P C^\infty(M)$.  Let $[u_k(g)]$ denote the $k$-th eigenspace of $\Delta_g$, viewed as a point of $\mathbb P C^\infty(M)$ whenever it is one-dimensional. For $D\subset\mathbb P C^\infty(M)$, define
\[
\mathcal G_+(k,D):=\{g\in\mathcal G_+:\lambda_k(g)\text{ is simple and }[u_k(g)]\in D\}.
\]

\begin{theoremB}\label{theorem b}
Let $k\geq1$ be a positive integer and let $D\subset\mathbb P C^\infty(M)$. If $D$ is open, dense, or residual, then $\mathcal G_+(k,D)$ is respectively open, dense, or residual in $\mathcal G_+$. In particular, if $D$ is residual, then the set of metrics $g$ such that every nonconstant eigenspace of $\Delta_g$ is simple and belongs to $D$ is residual in $\mathcal G_+$.
\end{theoremB}

Note that, this extends Uhlenbeck's theorem from the specific properties of being Morse and transverse to zero to \emph{arbitrary} generic properties of \emph{smooth} functions. Once \cref{theorem A} is known, the passage from generic properties of functions to generic properties of eigenfunctions is essentially topological: the preimage of a dense set under an open map is dense, and residual sets are handled by taking countable intersections.

Although our result is close in spirit to Uhlenbeck's, the proof uses a different mechanism. The argument in \cite{Uhlenbeck} is carried out in Banach spaces and relies on Sard--Smale transversality. In practice, this is particularly well suited to properties described by finite-order conditions. Our approach instead uses the Nash--Moser inverse function theorem in the category of tame Fr\'echet manifolds, following Hamilton\cite{Hamilton}. The analytic work is concentrated in proving a single openness theorem for the eigenfunction projection; once this is done, arbitrary residual subsets of $\mathbb P C^\infty(M)$ are treated simultaneously, including properties which are not naturally encoded by finitely many jets.
The proof of the analytic theorem has two local parts. Near a nondegenerate critical point, a finite-dimensional polynomial calculation solves the prescribed finite jet of the equation. In dimensions at least three, the remaining datum is solved by a characteristic integral along an attracting vector field. Finite-order vanishing and Lyapunov estimates give a fixed loss of derivatives and smooth tame dependence. In dimension two, the additional condition in $\mathcal N$ instead gives an explicit conformal solution.

Away from the critical points, a conformal variation corrects the mean of the datum. A traceless variation then reduces the remaining equation to a divergence equation, solved by a fixed support-preserving Bogovskii operator. Patching the two constructions gives the required smooth tame right inverse.

Uhlenbeck's theorem provides an open dense locus of simple eigenpairs where this analytic construction applies. The Nash--Moser theorem then proves openness of the eigenfunction projection. All subsequent genericity statements follow from this openness and elementary category arguments.

We conclude by describing the organization of the paper. In Section~1 we recall the analytic and geometric background used throughout: the Fr\'echet topologies on spaces of smooth functions and metrics, the basic spectral theory of the Laplace--Beltrami operator, and the Nash--Moser inverse function theorem in Hamilton's tame category.

Section~2 is devoted to the submersivity of the Laplace--Beltrami operator with respect to the metric. We first solve the linearized equation on jets at the critical points and prove the surjectivity of the homogeneous polynomial model. We then solve the residual equation near critical points by the characteristic method and Lyapunov estimates, and solve the equation away from critical points using the Bogovskii operator. Finally, we prove that the set $\mathcal N$ is open and dense. This is done by perturbing the second jets at the finitely many Morse critical points. We also show directly that the analogous surjectivity statement fails in dimension one.

In Section~3 we turn to eigenfunctions. We recall the part of Uhlenbeck's work needed in the sequel and prove that the simple eigenpair space is a tame Fr\'echet manifold. The submersivity theorem from Section~2 then gives Theorem~A. We use it to prove Theorem~B and its simultaneous version for all eigenvalues. We also discuss two extensions. First, restricting to metrics induced by embeddings into Euclidean space and using Nash's theory of free embeddings~\cite{Nash}, we obtain the analogous genericity statements for induced metrics. Second, we treat Schr\"odinger operators of the form
\[
\Delta_g+V,\qquad(g,V)\in\mathcal G_+\times C^\infty(M),
\]
proving an analogue of the flexibility theorem on an open dense subset of the corresponding eigenpair space and deriving the same genericity consequences. The appendix contains the coordinate computation of the metric differential and the Lyapunov estimates used in the characteristic construction.

We expect that the methods developed here can also be adapted to manifolds with boundary, after imposing suitable boundary conditions. We do not pursue this direction in the present paper.
\tableofcontents
\section{Preliminaries}
\subsection{Spaces of smooth functions and  metrics}
Throughout the paper, we fix a compact connected smooth manifold $M$ without boundary of dimension $n\ge 2$. Connectedness is retained throughout: the openness and Morse-density statements below do not extend verbatim to disconnected manifolds. The space of smooth real--valued functions on $M$ is denoted by $C^\infty(M)$. Fix a finite atlas $\{(U_i,\varphi_i)\}_{i=1}^N$ of $M$ and let $V_i\subset U_i$ be an open set with compact closure contained in $U_i$ and such that $\{V_i\}_{i=1}^N$ covers $M$. For each $k\ge 0$ we define the $C^k$ norm on $M$ as
\begin{equation}\label{C^k norm}
\|u\|_{C^k(M)} \coloneqq \max_{i}
\big\| u\circ \varphi_i^{-1}\big\|_{C^k(\varphi_i(V_i))},
\end{equation}
where on each $\varphi_i(V_i)\subset \mathbb{R}^n$ we use the standard $C^k$ norm. To simplify notation, we will often simply write $\|u\|_k=\|u\|_{C^k(M)}$ when there is no confusion. Observe that different choices of atlas give rise to equivalent norms. This is because the derivatives of the transition functions are bounded on a compact set. We endow  $C^\infty(M)$ with the topology generated by the countable family of norms $(\|\cdot\|_{C^k})_{k\in \mathbb{N}}$. With this topology, $C^\infty(M)$ is a Fr\'echet space; see \cite{Hamilton} for further details.
It will also be useful to consider Sobolev norms $\|\cdot\|_{H_k(M)}$ on $C^\infty(M)$, since these are well suited to regularity estimates; see \cref{schauder estimates}. Locally, for $u\colon U\subset \R^m\to \R$, the Sobolev norm is defined as
\[\|u\|_{H_k(U)}\coloneqq\left(\sum_{|\beta|\le k}\int_U\partial^\beta u(x)^2dx\right)^{\frac12}.\]
Using the same finite atlas, one defines the corresponding Sobolev norms on $M$. Despite not being equivalent, the families of norms $\{\|\cdot\|_{C^k(M)}\}_{k\ge 0}$ and $\{\|\cdot\|_{H_k(M)}\}_{k\ge 0}$ induce the same topology on $C^\infty(M)$. In particular, the Sobolev norms also define a Fr\'echet topology on $C^\infty(M)$. We will also denote by $H_k(M)$ the space of square--integrable functions with distributional derivatives up to order $k$. Observe that on a compact manifold being square integrable does not depend on the choice of the underlying metric.

We denote by $T^*M$ the cotangent bundle of $M$ and by $S^2T^*M$ the vector bundle of symmetric $(0,2)$--tensors on $M$. Among them, we denote by $S_+^2T^*M$ the subspace of the positive definite ones. We will be especially interested in sections of these space. Thus, we denote by $\cG$ the space of smooth sections of $S^2T^*M$ and by $\cG_+$ the space of smooth sections of $S^2_+T^*M$, which coincides with the space of smooth Riemannian metrics. Observe that $\cG_+$ is an open cone in $\cG$. As we did for smooth functions, we can define a Fr\'echet space structure on $\cG$ whose topology is generated by the $C^k$ norms. For each $g\in \cG_+$ we denote by $\mathrm{d}\vol_g$ the Riemannian volume measure associated with $g$ and by
\[\langle u_1,u_2\rangle_{L^2_g}
    \coloneqq \int_M u_1u_2 \,\mathrm{d}\vol_g, 
    \qquad
    \|u\|_{L^2_g}
    \coloneqq \bigl(\langle u,u\rangle_{L^2_g}\bigr)^{1/2},\]
the corresponding $L^2$–inner product and norm on $C^\infty(M)$. Observe that we could also define the norms $\|\cdot\|_{C^k_g(M)}$ using the metric $g$. Instead, we prefer to work with the (equivalent) fixed norm \eqref{C^k norm}, which is independent of $g$ and thus allows us to vary the metric while keeping the underlying Fr\'echet spaces and norms fixed. Moreover, we'll denote by $\cG_+^k$ the space of $C^k$ Riemannian metrics on $M$.

Analogously, we define a Fr\'echet space structure on the set of smooth vector fields on $M$ whose topology is generated by the $C^k$ norms.
\subsection{The Laplace--Beltrami operator}
Let  $g$ be a smooth Riemannian metric on $M$. We use the nonpositive sign convention for the Laplace--Beltrami operator:
\[
    \Delta_g \coloneqq \dive_g \nabla_g,
\]
so that $\Delta_g$ is nonpositive and self--adjoint on $H_2(M)$ with respect to the $L^2_g$ inner product on $L^2_g(M)$. Since $M$ is compact the spectrum of $\Delta_g$ consists of a  discrete sequence,
\[
0=-\lambda_0(g)>-\lambda_1(g)\ge -\lambda_2(g)\ge \cdots \to -\infty,
\]
where each eigenvalue is repeated according to its multiplicity. Thus, if $-\lambda$ is an eigenvalue of $\Delta_g$, the corresponding eigenfunctions are precisely the nonzero solutions of
\[
(\Delta_g+\lambda)u=0.
\] 
The eigenspace corresponding to $0$ consists of the constant functions. We shall mainly be interested in the negative eigenvalues. The multiplicity of an eigenvalue $-\lambda$  is denoted by $m_\lambda(g)\coloneqq \dim \ker(\Delta_g+\lambda)$ with the convention that $m_\lambda(g)=0$ if $-\lambda$ is not an eigenvalue of $\Delta_g$.

The range of the operator $\Delta_g+\lambda$ is described by the Fredholm alternative for selfadjoint operator. In particular if $g\in \cG_+$, then for every $\lambda\in\mathbb R$, the kernel ${\ker(\Delta_g+\lambda)}$ is finite-dimensional and
\[
(\Delta_g+\lambda)(C^\infty(M))
=
\ker(\Delta_g+\lambda)^{\perp_{L^2_g}}\cap C^\infty(M).
\] 
We shall use the following elliptic estimates. 

\begin{theorem}[{\cite[Theorem 2, Chapter 6.3]{Evans}}]\label{schauder estimates}
    Let $M$ be a compact manifold and $L$ an elliptic operator with $C^{k+1}$ coefficients. Fix a scalar $\lambda$, a positive integer $k\geq 0$. There exists a positive constant $C=C(L)$ depending continuously on $L$ such that
    \[\|u\|_{H_{k+2}(M)}\le C\left(\|Lu\|_{H_{k}(M)}+\|u\|_{L^2(M)}\right),\]
     for every $u\in H_{k+2}(M)$. Moreover if $g$ is smooth, $u\in H_1(M)$ and $Lu\in C^\infty(M)$, then have $u\in C^\infty(M)$.
\end{theorem}
\begin{remark}
    \cref{schauder estimates} holds for operators of the Laplace--Beltrami type. However, not that if $g$ is a $C^k$ metric, then $\Delta_g$ has $C^{k-1}$ coefficients. This is an immediate consequence of tis definition.
\end{remark}
Throughout this paper, we will always assume that $k$ is large enough that $H_k(M)$ embedds into $C^2(M)$. Such an embedding is given by \cite[Theorem 6, Chapter 5.6.3]{Evans} for $k$ greater than a fixed constant depending only on the dimension of $M$.

\subsection{Tame Fr\'echet geometry and the Nash--Moser theorem.}\label{tame geometry}
In this subsection we recall the notions of tame Fr\'echet geometry needed in the sequel. The purpose is twofold: first, to fix the notions of tame maps and tame submanifolds used below; second, to state the Nash--Moser inverse function theorem in the form in which it will be applied.

We begin with the underlying Fr\'echet-manifold structure, see \cite[Definition 4.1.1 Sec 1]{Hamilton}. 
\begin{definition}[Fr\'echet manifold]
    Let $E$ be a Fr\'echet space. A \emph{Fr\'echet manifold} modeled on $E$ is a Hausdorff topological space endowed with an atlas of charts taking values in open subsets of $E$, whose transition maps are smooth in the Fr\'echet sense. 
\end{definition}
For the Nash--Moser inverse function theorem, however, this ordinary Fr\'echet structure is not sufficient. One must work with Fr\'echet spaces and maps satisfying an additional tameness assumptions. We will use the framework of tame Fr\'echet spaces in the sense of Hamilton, see \cite{Hamilton}, recalling only the parts of the theory that will be used explicitly.
In particular, all Fr\'echet spaces considered below are assumed to be graded: their topology is induced by a family of increasing seminorms $\|\cdot\|_k$, satisfying
\[
\|f\|_k \leq \|f\|_{k+1}
\]
for every $f$ and every $k$. The additional condition that a graded Fr\'echet space be tame is a compatibility condition with a model space of rapidly decreasing sequences. Refer \cite[Definition 1.3.2 Sec 2]{Hamilton} for the precise definition. The standard examples relevant here, such as spaces of smooth sections over compact manifolds endowed with their usual $C^k$-norms, are tame Fr\'echet spaces.

The notion of tameness for maps is expressed by estimates with a controlled loss of derivatives. Since we will have to verify such estimates later, we now recall the relevant definitions from \cite{Hamilton}.
\begin{definition}[Tame maps]
Let $F$ and $G$ be Fr\'echet spaces, and let $L\colon F \to G$ be a linear map. We say that $L$ is a \emph{tame linear map} if there exist $r,b\in \mathbb N$ such that
\[
\|Lf\|_k\leq C_k\,\|f\|_{k+r}
\]
for every $k\geq b$, where the positive constant $C_k$ depends only on $k$.
\\Let $U\subseteq F$ be open, and let $P\colon U \to G$ be a nonlinear continuous map. We say that $P$ satisfies a \emph{tame estimate} if there exist $r,b\in \mathbb N$ such that
\[
\|P(f)\|_k \leq C_k\bigl(1+\|f\|_{k+r}\bigr)
\]
for every $k\geq b$, where the constant $C_k$ depends only on $k$. The map $P$ is said to be \emph{tame} if it satisfies a tame estimate in a neighborhood of each point. A smooth tame map is a smooth map such that all its derivatives are tame. We say that the image of a smooth tame embedding is a smooth tame submanifold.
\end{definition}
\begin{remark}
    In practice, we will use the following simple criterion. If $X$ is a smooth tame submanifold of a tame Fr\'echet space $E$, then to prove that a map $P\colon X\to F$ is tame it is enough to prove, locally on $X$, estimates of the form
\[
\|P(f)\|_k \leq C_k(1+\|f\|_{k+r}),
\]
where the seminorms on $X$ are induced from the ambient space $E$.
\end{remark}
All partial differential operator are smooth and tame, in particular
the Laplace--Beltrami operator
\begin{align*}
    \Delta\colon \cG_+\times \R\times C^\infty(M)\to& C^\infty(M)
    \\
    g,\lambda,u\mapsto &(\Delta_g+\lambda) u
\end{align*}
is smooth and tame. This is an immediate consequence of \cite[Corollary 2.2.7 Sec 2]{Hamilton}.

The verification of tame estimates will often rely on standard interpolation inequalities. These inequalities allow us to estimate intermediate norms in terms of lower and higher norms, and to control products without losing more derivatives than necessary. We record the precise form that will be used below.
\begin{proposition}[{\cite[Theorem 2.2.1, Corollary 2.2.3 Sec 2]{Hamilton}}]\label{interpolation}
    Let $M$ be a compact smooth manifold and $f,g\in C^{\infty}(M)$. If $k\leq \ell\leq m$, then
    \[
    \|f\|_\ell^{m-k}\leq C_{k,\ell,m}\|f\|_k^{m-\ell}\|f\|_m^{\ell-k},\]
    where $C_{k,\ell,m}$ is a positive constant depending on $M$, $k,\ell$ and $m$. This estimate also holds for Sobolev norms, possibly with different constants, and for smooth vector fields. Moreover, we have 
    \[
    \|fg\|_m\leq C_{0,m}(\|f\|_0\|g\|_m+\|f\|_m\|g\|_0),\]
    where $C_{0,m}$ is a positive constant depending on $M$ and $m$. 
\end{proposition}

We will need to show that families of inverse maps are smooth and tame. The following criterion reduces this to proving continuity and tameness for the inverse family.
\begin{proposition}[{\cite[Theorem 3.1.1 Sec 3]{Hamilton}}]\label{conti smooth}
    Let $L:(U\subseteq E)\times F\to G$ be a smooth tame family of linear maps. Suppose that the equation $L(e)f = g$ has a unique solution $f$ for all $e$ and $g$. Assume also that the family of inverses $V(e)g = f$ is continuous and tame as a map $V\colon U\times G \to F$. Then $V$ is a smooth tame map.
\end{proposition}
We are now ready to state the Nash--Moser inverse function theorem. The main difference from the classical Banach implicit function theorem, besides the additional tameness requirement, is that one must assume the differential to be invertible on a whole neighborhood of the point, rather than only at the point itself.

\begin{theorem}[{\cite[Theorem 1.1.1, 1.1.2, 1.1.3 Sec 3]{Hamilton}}]\label{nash moser}
    Let $F,G$ be tame Fr\'echet spaces. Let $U$ be an open subset of $F$ and $P:U\subseteq F \to G$ be a smooth tame map. Then
    \begin{enumerate}
        \item Assume that $DP\colon U\times F \to G$ has a smooth tame right inverse, i.e. $RP\colon U\times G \to F$ such that for each $p\in U$, $f\in F$ we have $D_pP(R_pP(f))=f$. Then $P$ is locally surjective, i.e. in a neighborhood of any point it has a smooth tame right inverse. 
        \item Assume that $DP\colon U\times F\to G$ has a smooth tame left inverse. Then $P$ is locally injective with a smooth tame left inverse. 
        \item Assume that $DP\colon U\times F\to G$ has a  smooth tame inverse. Then $P$ is locally invertible and the local inverse is a smooth tame map.
    \end{enumerate} 
\end{theorem}
\section{Submersivity of the Laplacian}\label{section N}
We can interpret the Laplacian operator as a smooth map
\[\Delta:\cG_+\times C^\infty(M)\to C^\infty(M),\]
where $\Delta(g,u)=\Delta_g(u)$. In this section, we'll study the differential of the Laplacian with respect to the metric component. To simplify notation, we'll denote such differential as $\dot\Delta_g[\cdot](u)$ at a point $(g,u)$. We first introduce the space we'll be working with.
\begin{definition}[Submersivity set]\label{nonresonant set}
    Let $M$ be a compact smooth manifold of $\dim(M)\geq 2$ and let $k$ be a sufficiently large positive integer. We define the following subset of $\cG_+^{k}\times H_{k}(M)$ as
    \[\mathcal N_{C^k}=\begin{cases}
    \{(g,u):u\ \mathrm{\ Morse},\ \mathrm{crit}(u)\cap(\Delta_gu)^{-1}(0)=\emptyset\}\quad\mathrm{if\ }\dim(M)=2\\
    \{(g,u):u\text{ Morse}\}\quad\quad\quad\quad \quad\quad\ \quad\quad\quad\quad\ \quad\quad\mathrm{if\ }\dim(M)>2.
\end{cases}\]
We will denote its smooth part by
\[\mathcal N=\mathcal N_{C^k}\cap(\cG_+\times C^\infty(M))\subseteq \cG_+\times C^\infty(M).\]
\end{definition}
The aim of this section is to prove the following theorem, which will serve as the main tool prove the results in the next section.
\begin{theorem}\label{laplacian use}
    Let $M$ be a smooth compact connected manifold with $\dim(M)\geq 2$ and let $(g_0,u_0)\in\mathcal N$. Then there exists a neighborhood $\mathcal{V}$ of $(g_0,u_0)$ in $\cG_+\times C^\infty(M)$ such that the map
    \[\dot{\Delta}\colon\mathcal V\times \cG\to C^{\infty}(M)\] is surjective with a smooth tame right inverse. Moreover, the set $\mathcal N$ is open dense. 
\end{theorem}
 As an immediate consequence, we obtain the following corollary.
\begin{corollary}\label{dimension 3 cor}
    Let $M$ be smooth compact connected manifold with $\dim(M)\geq 3$. Fix a Morse function $u\in C^\infty(M)$. Then the map
    \[\Delta(u):\cG_+\to C^{\infty}(M)\]
    defined as $\Delta(u)(g)=\Delta_g(u)$ is open.
\end{corollary}
\begin{proof}
    This result is an immediate consequence of \cref{nash moser} and \cref{laplacian use}.
\end{proof}
Proving \cref{laplacian use} amounts to study the equation
\begin{equation}\label{diff equation}
    \dot\Delta_g[\cdot](u)=v
\end{equation}
when $u$ is Morse and $v$ is any smooth function. To do so, it will be convenient to give an explicit expression for the differential of the Laplacian. We will often omit the Riemannian metrics when denoting gradients if there is no ambiguity.

Recall that a smooth metric $g$ on $M$ induces a fiber metric for all tensor bundle, see \cite[Lemma 3.1]{MR1468735}. In particular for two smooth $(0,2)$--tensors  $h_1,h_2$ we have $$\langle h_1,h_2\rangle_g=g^{ik}g^{j\ell}(h_1)_{ij}(h_2)_{k\ell}.$$ Also recall that one can define the divergence of a $(0,2)$--tensor $h$ using a metric $g$. It is given by the $(0,1)$--tensor $\mathrm{div}_g(h)=\mathrm{tr}_g(\nabla h)$. Moreover, if $\alpha$ is a $(0,1)$--tensor and $g$ is a Riemannian metric, we denote by $\alpha^\sharp$ the unique vector field such that
$\langle\alpha^\sharp,\cdot\rangle_g=\alpha$.
\begin{lemma}[Expression for the differential]\label{lem:metric variation of laplacian}
    Let $(M,g)$ be a Riemannian manifold. Then
    \begin{equation}\label{differential of Lapl}
    \dot{\Delta}_g[h](u)=-\langle h,\nabla^2 u\rangle_g+\langle\frac12\nabla(\mathrm{tr_g(h)})-\mathrm{div}_g(h)^\sharp,\nabla u\rangle_g.\end{equation}
    In particular, $\dot\Delta$ is smooth tame in $q$ and $h$. For a conformal variation $h_{ij}=\psi g_{ij}$, then \begin{equation}\label{conformal equation}
    \dot{\Delta}_g[\psi g](u)=-\psi \Delta_g u+ \frac{n-2}2\langle \nabla \psi,\nabla u\rangle_g.\end{equation}
    Moreover, for a traceless variation $g^{ij}h_{ij}=0$ we have
    \begin{equation}\label{traceless equation}
    \dot{\Delta}_g[h](u)=-\dive_g( h(\nabla u,\cdot))
    \end{equation}
    Finally, let $W$ be a coordinate chart and $\{dx^i\}_{i=1}^n$ be the corresponding coordinate $1$--forms. Let $F=F_{ij}dx^i\otimes dx^j$ be a symmetric $(0,2)$--tensor whose entries are constant in coordinates and $w$ be a function compactly supported in $W$. Then
    \begin{equation}\label{final 3.8}
        \dot{\Delta}_g[wF](u)=a_qw+b_q(w),
    \end{equation}
    where
    \[a_q=-\langle F,\nabla^2 u \rangle_g+\langle\frac 12\nabla\tr_g F-\dive_g F^\sharp,\nabla u\rangle_g,\quad b_q=\frac12\tr_g F\nabla u-F(\nabla u,\cdot)^\sharp.\]
    Moreover, $a_q$ and $b_q$ are smooth tame in $q$.
\end{lemma}
\begin{proof}
    See the appendix.
\end{proof}
\begin{remark}[Conformal equation in dimension $2$]
    Observe that \eqref{conformal equation} becomes significantly simpler when $n=2$. This will be used in the proof of \cref{sol in crit} in dimension $2$.
\end{remark}Observe that \cref{laplacian use} does not hold in dimension $1$. Indeed, the differential of the Laplace--Beltrami operator is never surjective in this case.
\begin{theorem}\label{no 1}
    Let $M$ be a smooth compact connected manifold with $\dim(M)=1$. Given a nonconstant function $u\in C^\infty(M)$, the metric differential $\Delta_g[\cdot](u)$ is not surjective in $(g,u)$ for every $g\in\cG_+$.
\end{theorem}
\begin{proof}
    Since every connected smooth compact manifold is diffeomorphic to the circle, we can assume $M=S^1$. Fix $(g,u)\in\cG_+\times C^\infty(S^1)$. Without loss of generality, we can assume $\vol_g(S^1)=1$. Moreover, by fixing arclength coordinates for $g$, we can assume that $g$ is the standard Euclidean metric. Also assume that we have fixed an orientation of $S^1$, so that it makes sense to treat $g$ as a function on $S^1$. Since there are no nonzero traceless metric variations in dimension $1$, the expression for $\Delta_g[\cdot](u)$ reduces to the conformal one
    \[\Delta_g[\psi g](u)=-\psi \Delta_g u-\frac 12\langle \nabla \psi,\nabla u\rangle_g.\]
    Under our assumptions, it further simplifies to
    \[\Delta_g[\psi g](u)=-\psi u''-\frac 12\psi'u'=-\frac{((u')^2\psi)'}{2u'}.\]
    We want to show that there exists $f\in C^\infty(S^1)$ such that
    \[((u')^2\psi)'=2fu'\]
    has no solution $\psi$. Let $\Omega\subseteq S^1$ be the set of points where $u'$ is nonzero. Since $u$ is nonconstant, this is a nonempty open set. Fix a connected component $\tilde\Omega$ of $\Omega$. Since $u$ is not constant, it has at least $2$ critical points. In particular, the boundary of $\tilde\Omega$ is given by two distinct points $x_1$ and $x_2$. For every $\psi\in C^\infty(S^1)$ we have
    \[\int_{\tilde{\Omega}}2u'\Delta_g[\psi g](u)=-\int_{\tilde{\Omega}}((u')^2\psi)'=\pm(u')^2(x_2)\psi(x_2)\mp(u')^2(x_1)\psi(x_1)=0.\]
Thus, if $\Delta_g[\cdot](u)=f$ has a solution, then $f$ must verify the equation
\[\int_{\tilde{\Omega}}u'f=0.\]
This is clearly a nontrivial condition, thus $\Delta_g[\cdot](u)$ is not surjective.
\end{proof}
\begin{remark}
    A finite regularity version of \cref{no 1} can be proven in the same way.
\end{remark}
\subsection{Solution on jets of the linearized Laplace--Beltrami operator at critical points}
The behavior of \eqref{diff equation} is very different around critical points of $u$. Thus, we split the proof of \cref{laplacian use} in $2$ parts. First, we will build a smooth tame right inverse in a neighborhood of the critical points of $u_0$. In particular, we prove the following.
\begin{proposition}[Solution around critical points]\label{sol in crit}
    Let $q_0=(g_0,u_0)\in \mathcal N$. Then there exists an open set $W_{\crit}\subseteq M$ with smooth boundary, a neighborhood $\mathcal{V}$ of $q_0$ in $\cG_+\times C^\infty(M)$, and a smooth tame map $$R_{\crit}\colon\mathcal V\times C^\infty(M)\to \cG$$ as follows. For every $(g,u)\in\mathcal{V}$, the set $\crit (u)$ is contained in a fixed compact subset of $W_{\crit}$ with smooth boundary. Moreover,
    \begin{equation}\label{new to cref}
        \dot{\Delta}\circ R_{\crit}v(x)=v(x)\quad \text{for every $x\in W_{\crit},v\in C^\infty(M)$.}
    \end{equation}
\end{proposition}
We will prove \cref{sol in crit} by constructing an appropriate solution of \eqref{diff equation} using the method of characteristics. However, such a solution only converges when the $d$--th jet of the datum $v$ in the critical points of $u$ is $0$ for a fixed $d$ sufficiently large. In this subsection, we reduce the problem to this case by solving \eqref{diff equation} in its jet. We first record the smooth and tame continuity of the critical points of $u$.
\begin{lemma}[Persistence of Morse critical points]
\label{lem:persistence of critical points}
Let $u_0\in C^\infty(M)$ be Morse, with critical points
\[
\crit(u_0)=\{x_1,\ldots,x_m\}.
\]
Then there exist pairwise disjoint open neighborhoods $W_i\ni x_i$, a neighborhood $\mathcal V$ of $u_0$, and smooth tame maps
\[
x_i\colon\mathcal V\longrightarrow W_i, \qquad x_i(u_0)=x_i,\]
such that, for every $u\in\mathcal V$,
\[
\crit(u)=\{x_1(u),\ldots,x_m(u)\},
\]
and each $x_i(q)$ is a nondegenerate critical point of $u$. In particular, $x_i(q)$ is the unique critical point of $u$ in $W_i$.
\end{lemma}
\begin{proof}
Choose the neighborhoods $W_i$ pairwise disjoint and so small that $x_i$ is the only critical point of $u_0$ in $\overline{W_i}$. Working in a chart and writing the equation $Du(x)=0$ in coordinates, the Morse condition says that the differential of the left-hand side with respect to $x$, namely \ $\nabla^2u_0(x_i)$, is invertible. The implicit function theorem with parameters therefore produces, for $q$ in a neighborhood of $q_0$, a unique solution $x_i(q)\in W_i$, depending smoothly on $q$. Since $x_i(q)$ is obtained by solving a finite-dimensional equation whose coefficients are finitely many derivatives of $u$, the map $q\mapsto x_i(q)$ is smooth tame. Shrinking $\mathcal V$ further, $\nabla^2u(x_i(q))$ stays invertible by continuity.

Finally, $Du_0$ does not vanish on the compact set $M\setminus\bigcup_iW_i$; hence $Du$ does not vanish there either, provided $u$ is close enough to $u_0$ in $C^1(M)$.
\end{proof}
Let $u$ be any Morse function and $x_u$ be a nondegenerate critical point of $u$. Choose geodesic normal coordinates centered at $p$. For simplicity, identify $p$ with the origin of $\R^n$. Set
\[
H\coloneqq\nabla^2u(0),
\]
which at a critical point coincides with the Euclidean Hessian of the coordinate representation of $u$. Since $u$ is Morse, it is a nonsingular symmetric matrix. In these coordinates,
\[
g_{ij}(x)=\delta_{ij}+O(|x|),
\qquad
Du(x)=Hx+O(|x|^2),
\qquad
\nabla^2u(x)=H+O(|x|).
\]

For $d\geq0$, let $\mathcal P_d$ denote the space of homogeneous polynomials of degree $d$ on $\R^n$. Let $\mathcal S_d$ denote the space of symmetric-matrix-valued homogeneous polynomials of degree $d$, that is, of tensors $h=(h_{ij})$ with $h_{ij}=h_{ji}\in\mathcal P_d$. If $h\in\mathcal S_d$, then
\begin{equation}
\operatorname{tr}_gh=\operatorname{tr}h+O(|x|^{d+1}),\qquad\dive_gh=\dive h+O(|x|^{d}),
\end{equation}
where the operations on the right-hand side are the Euclidean ones. Substituting these expansions into \eqref{differential of Lapl} gives
\begin{equation}\label{expansion of the linearization}
\dot\Delta_g[h](u)= A_{H,d}h+O(|x|^{d+1}), \qquad h\in\mathcal S_d,
\end{equation}
where
\begin{equation}\label{definition of A_Hd}
 A_{H,d}h\coloneqq-\langle H,h\rangle+\tfrac12\bigl\langle Hx,D(\operatorname{tr}h)\bigr\rangle - \bigl\langle\dive h,Hx\bigr\rangle .
\end{equation}
This is because $D(\operatorname{tr}h)$ and $\dive h$ are homogeneous of degree $d-1$, while $Hx$ is homogeneous of degree $1$. Thus, $ A_{H,d}h$ is exactly the homogeneous component of degree $d$ of
$\dot\Delta_g[h](u)$, and \eqref{definition of A_Hd} defines a linear map
\[
 A_{H,d}\colon\mathcal S_d\longrightarrow\mathcal P_d .
\]
Observe that $ A_{H,d}$ depends only on $H$: all the information carried by $g$ and by the higher-order Taylor coefficients of $u$ has been discarded. This is the model to which the equation reduces at the critical point, and the whole subsection rests on the following algebraic fact.
\begin{proposition}[The homogeneous polynomial model]
\label{homogeneous polynomial model}
Let $H\in\sym(n)$ be nonsingular, $n\geq 2$. Then
\[
 A_{H,d}\colon\mathcal S_d\longrightarrow\mathcal P_d
\]
is surjective for every $d\geq0$. Moreover, if $H_0$ is fixed, right inverses $R_{H,d}$ of $ A_{H,d}$ can be chosen so as to depend smoothly on $H$ in a neighborhood of $H_0$, with locally uniformly bounded operator norm.
\end{proposition}

\begin{proof}
For $Q\in O(n)$, define $(Q^*h)(x)=Q^Th(Qx)Q$. A direct change of coordinates gives
\[
A_{Q^THQ,d}(Q^*h)=(A_{H,d}h)\circ Q.
\]
We may therefore assume that
\[
H=\operatorname{diag}(\lambda_1,\ldots,\lambda_n),\qquad\lambda_i\neq0\qquad\text{for every }i .
\]
Since the monomials $x^\alpha$ with $|\alpha|=d$ span $\mathcal P_d$, it suffices to prove that each of them lies in the image of $ A_{H,d}$. Throughout the proof we set
\[
S_\alpha\coloneqq\sum_{k=1}^n\lambda_k\alpha_k .
\]

We first consider diagonal tensors. Fix $i$ and let $h\in\mathcal S_d$ be the tensor whose only nonvanishing component is $h_{ii}=x^\alpha$. Then
\[
\operatorname{tr}h=x^\alpha,
\qquad
\langle H,h\rangle=\lambda_ix^\alpha,
\qquad
(\dive h)_j=\delta_{ij}\,\partial_ix^\alpha ,
\]
hence
\[
\tfrac12\langle Hx,Dx^\alpha\rangle=\tfrac12S_\alpha\,x^\alpha,
\qquad
\langle\dive h,Hx\rangle=\lambda_i\alpha_i\,x^\alpha .
\]
Therefore
\begin{equation}\label{diagonal coefficient}
 A_{H,d}h=c_i(\alpha)\,x^\alpha,
\end{equation}
where $c_i(\alpha)\coloneqq\tfrac12S_\alpha-\lambda_i(\alpha_i+1)$. If $c_i(\alpha)\neq0$ for some $i$, then $x^\alpha$ belongs
to the image of $ A_{H,d}$. Call $\alpha$ \emph{exceptional} for $H$ when this fails, that is when
\begin{equation}\label{exceptional monomial}
\lambda_i(\alpha_i+1)=\tfrac12S_\alpha
\qquad\text{for every }i=1,\ldots,n .
\end{equation}

Assume now that $\alpha$ is exceptional. We want to show that $x^\alpha$ belongs to the image of $ A_{H,d}$. Since the $\lambda_i$ are nonzero, \eqref{exceptional monomial} implies $S_\alpha\neq0$. Multiplying the $i$-th identity by
$\frac{\alpha_i}{\alpha_i+1}$ gives
\[
\lambda_i\alpha_i=\frac{S_\alpha}{2}\,\frac{\alpha_i}{\alpha_i+1}.
\]
Summing over $i$ and dividing by $S_\alpha\neq 0$ gives
\begin{equation}\label{exceptional identity}
\sum_{i=1}^n\frac{\alpha_i}{\alpha_i+1}=2 .
\end{equation}
Each summand lies in $[\tfrac12,1)$ when $\alpha_i\geq1$ and vanishes otherwise, so \eqref{exceptional identity} has at least three and at most four nonzero entries. In particular $n\geq3$, $|\alpha|\geq 4$, and there is an index $j$ with $\alpha_j\geq1$.

Fix a $j$ as above, and pick any index $i\neq j$. Define
\[
\gamma\coloneqq\alpha+e_i-e_j\in\mathbb N^n,
\]
and let $h\in\mathcal S_d$ be the tensor whose only nonvanishing components are $h_{ij}=h_{ji}=x^\gamma$. Since $H$ is diagonal and $h$ is
off-diagonal,
\[
\operatorname{tr}h=0,\qquad\langle H,h\rangle=0,
\]
while
\[
(\dive h)_i=\partial_jx^\gamma,\qquad
(\dive h)_j=\partial_ix^\gamma,\qquad
(\dive h)_k=0\ \ (k\neq i,j).
\]
If $\alpha_j=1$, direct computation gives $A_{H,d}h=-\lambda_j(\alpha_i+1)x^\alpha$, which proves the claim. Thus, suppose $\alpha_j\ge2$ and define $\beta=\alpha+2e_i-2e_j\in\mathbb N^n$. Using $x_i\,\partial_jx^\gamma=(\alpha_j-1)x^{\beta}$ and
$x_j\,\partial_ix^\gamma=(\alpha_i+1)x^{\alpha}$, we get
\begin{equation}\label{off diagonal image}
 A_{H,d}h=-\lambda_j(\alpha_i+1)\,x^\alpha-\lambda_i(\alpha_j-1)\,x^{\beta}.
\end{equation}
The coefficient $a=\lambda_j(\alpha_i+1)$ is nonzero, and $|\beta|=d$.

We claim that $\beta$ is not exceptional. Suppose it were. Since $n\geq3$ we may choose $k\notin\{i,j\}$, which gives $\beta_k=\alpha_k$. Then \eqref{exceptional monomial} applied to $\beta$ and to $\alpha$ gives
\[
\tfrac12S_\beta=\lambda_k(\beta_k+1)=\lambda_k(\alpha_k+1)=\tfrac12S_\alpha .
\]
On the other hand, the identities for the index $i$ give
\[
\tfrac12S_\beta=\lambda_i(\beta_i+1)=\lambda_i(\alpha_i+3),
\qquad
\tfrac12S_\alpha=\lambda_i(\alpha_i+1).
\]
This implies $2\lambda_i=0$, a contradiction since $H_0$ is nonsingular. Therefore $c_\ell(\beta)\neq0$ for some $\ell$, and Step 1 provides $\widetilde h\in\mathcal S_d$ with $ A_{H,d}\widetilde h=x^\beta$. Setting $b\coloneqq\lambda_i(\alpha_j-1)$, \eqref{off diagonal image} yields
\[
 A_{H,d}\Bigl(-\tfrac1a\bigl(h+b\,\widetilde h\bigr)\Bigr)=x^\alpha .
\]
This proves that every monomial of degree $d$ is in the image, and hence that $ A_{H,d}$ is surjective.

We now prove the last part of the statement. Fix a right inverse $R_{H_0,d}$ of $ A_{H_0,d}$. The map $H\mapsto A_{H,d}$ is linear, hence continuous, so for $H$ close to $H_0$ the endomorphism $ A_{H,d}R_{H_0,d}$ of $\mathcal P_d$ is close to the identity and therefore invertible. Consequently
\[
R_{H,d}\coloneqq R_{H_0,d}\bigl( A_{H,d}R_{H_0,d}\bigr)^{-1}
\]
is a right inverse of $ A_{H,d}$; it depends smoothly on $H$ and its operator norm is bounded on a neighborhood of $H_0$.
\end{proof}
We now pass from homogeneous polynomials to jets. For a smooth function or tensor field $F$ defined near the origin, we denote by $[F]_d$ the homogeneous component of degree $d$ of its Taylor expansion at the origin, and by $j^N_xF$ the $N$-jet of $F$ at $x$. Let $\mathcal P_{\leq N}$ and $\mathcal S_{\leq N}$ denote respectively scalar and symmetric-tensor-valued polynomials of degree at most $N$. From now till the end of this subsection, we'll implicitly always used the coordinates above.

\begin{lemma}[Triangularity of the linearization]
\label{lem:triangularity}
Let $u\in C^\infty(M)$ and $x_u$ be a nondegenerate critical point of $u$. If $h=h_0+\cdots+h_N$ with $h_d\in\mathcal S_d$, then for every $0\leq d\leq N$
\[
\bigl[\dot\Delta_g[h](u)\bigr]_d= A_{H,d}h_d+\Phi_d(h_0,\ldots,h_{d-1}),
\]
where $\Phi_d$ is linear in its arguments and its coefficients are polynomial expressions in the Taylor coefficients of $g$, $g^{-1}$ and $u$ at $p$, of order at most $d+2$.
\end{lemma}

\begin{proof}
Since $Du(p)=0$, the operator $h\mapsto\dot\Delta_g[h](u)$ does not decrease the order of vanishing at $p$: the zero-order term $-\langle h,\nabla^2u\rangle_g$ clearly does not, and in the first-order term the derivatives of $h$ are multiplied by $Du=O(|x|).$ Hence each $h_e$ contributes only to the homogeneous components of degree $\geq e$. By \eqref{expansion of the linearization} its contribution in degree exactly $e$ is $ A_{H,e}h_e$. Collecting the terms of degree $d$ gives the stated formula, with $\Phi_d$ accounting for the contributions of $h_0,\ldots,h_{d-1}$. The last assertion follows by inspection of \eqref{differential of Lapl}, since $\dot\Delta_g[h](u)$ involves $g$, $g^{-1}$, $Dg$, $Du$ and $\nabla^2u$.
\end{proof}

\begin{proposition}[Solution on jets]\label{jet solve}
Let $q_0=(g_0,u_0)\in\mathcal N$, and let $N\geq0$. Then there exist a neighborhood $\mathcal V$ of $q_0$ in $\mathcal N$ and a map
\[
R_{\mathrm{jet},N}\colon\mathcal V\times C^\infty(M)\longrightarrow\cG,
\]
linear in the second variable and smooth tame, such that
\begin{equation}\label{jet equation}
j^N_x\Bigl(\dot\Delta_g\bigl[R_{\mathrm{jet},N}(q,v)\bigr](u)\Bigr)=j^N_xv
\qquad\text{for every }x\in\crit(u),
\end{equation}
for every $q=(g,u)\in\mathcal V$ and every $v\in C^\infty(M)$. Moreover, $R_{\mathrm{jet},N}(q,v)$ is supported in a fixed compact neighborhood of $\crit(u_0)$, and for every $m\geq0$ there is a constant $C_{m,N}>0$, independent of $q\in\mathcal V$, such that
\begin{equation}\label{jet estimate}
\bigl\|R_{\mathrm{jet},N}(q,v)\bigr\|_{C^m(M)}
\leq
C_{m,N}\,\|v\|_{C^N(M)} .
\end{equation}
\end{proposition}

\begin{proof}
We first work near a single critical point. Let $W_i$, $\mathcal V$ and $x_i\colon\mathcal V\to W_i$ be as in \cref{lem:persistence of critical points}, and fix $i$. Let $x$ be a fixed system of $g_0$-normal coordinates centered at $x_i(q_0)$ and defined on $W_i$, and for $q\in\mathcal V$ use the translated coordinates
\[
y\coloneqq x-x\bigl(x_i(q)\bigr),
\]
which are centered at the critical point $x_i(q)$ of $u$. Define
\[
 L_N(q)\colon\mathcal S_{\leq N}\longrightarrow\mathcal P_{\leq N},
\]
 where $L_N(q)h$ is the Taylor polynomial of degree $N$ of $\dot\Delta_g[h](u)$ at $y=0$.

Write $h=h_0+\cdots+h_N$ and $P=P_0+\cdots+P_N\in\mathcal P_{\leq N}$. By \cref{lem:triangularity}, the equation $ L_N(q_0)h=P$ is equivalent to the triangular system
\[
 A_{H,d}h_d=P_d-\Phi_d(h_0,\ldots,h_{d-1}),
\qquad
d=0,\ldots,N,
\]
where $H=\nabla^2u_0(x_i(q_0))$. Since $u_0$ is Morse, $H$ is nonsingular, and \cref{homogeneous polynomial model} solves each equation of the system. Solving them recursively, starting from $d=0$, produces $h$ with $ L_N(q_0)h=P$. Hence $ L_N(q_0)$ is surjective; let $R^0_N$ be a right inverse.

All the spaces involved are finite-dimensional, and by \cref{lem:triangularity} the operator $ L_N(q)$ depends only on the jets $j^{N+1}_{x_i(q)}g$ and $j^{N+2}_{x_i(q)}u$, evaluated at the point $x_i(q)$, which varies smoothly and tamely with $q$ by \cref{lem:persistence of critical points}. Consequently $q\mapsto L_N(q)$ is smooth, and continuous for the $C^{N+2}$ topology. Shrinking $\mathcal V$ to a $C^{N+2}$-ball around $q_0$, we may assume that $ L_N(q)R^0_N$ is invertible on $\mathcal P_{\leq N}$ for every $q\in\mathcal V$, with inverse bounded uniformly in $q$. Then
\[
R_N(q)\coloneqq R^0_N\bigl( L_N(q)R^0_N\bigr)^{-1}
\]
is a right inverse of $ L_N(q)$, depending smoothly on $q$, with
operator norm bounded uniformly on $\mathcal V$.

Let $\rho_i$ be a cutoff function supported in $W_i$ and identically $1$ on a neighborhood $K_i$ of $x_i(\mathcal V)$, fixed once and for all, and set
\[
R_{\mathrm{jet},N}(q,f)\coloneqq \sum_{i=1}^m\rho_i\cdot R_N^{(i)}(q)\Bigl(j^N_{x_i(q)}f\Bigr),
\]
where $R^{(i)}_N(q)$ is the operator constructed above at the $i$-th critical point, and where the polynomial tensor $R^{(i)}_N(q)(j^N_{x_i(q)}f)$, written in the coordinates $y$, is regarded as a smooth symmetric tensor on $W_i$ and extended by zero. The map is linear in $f$ by construction, and \eqref{jet equation} holds because $\rho_i$ is identically $1$ near $x_i(q)$ and the neighborhoods $W_i$ are pairwise disjoint.

The tensor $R_{\mathrm{jet},N}(q,f)$ is a polynomial of degree at most $N$ multiplied by a fixed cutoff, and its coefficients depend linearly on $j^N_{x_i(q)}f$, with a matrix bounded uniformly on $\mathcal V$. Since all the $C^m$ norms on the finite-dimensional space $\mathcal S_{\leq N}$ restricted to the fixed compact sets $\operatorname{supp}\rho_i$ are equivalent, and since $|j^N_xf|\leq\|f\|_{C^N(M)}$, we obtain \eqref{jet estimate}. In particular, $R_{\mathrm{jet},N}$ is a tame map, with loss of $N$ derivatives. Its derivatives with respect to $q$ are again of the same form, namely finite-dimensional expressions in the jets of $g$ and $u$ at $x_i(q)$ composed with the jet evaluation map; hence they satisfy analogous estimates, and $R_{\mathrm{jet},N}$ is smooth tame.
\end{proof}
\subsection{Right inverse of the linearized Laplace--Beltrami operator at critical points}
We can now give a complete solution of \eqref{diff equation} in a neighborhood of critical point of $(g_0,u_0)$. We will assume $\dim(M)>2$ for now. The $2$--dimensional case is significantly different and will be treated separately.

A key step in proving smoothness of the solution will be estimating appropriately the flow of a vector field. These estimates are given by Lyapunov theory, and are recorded in the following lemma.
\begin{lemma}[Lyapunov estimates]\label{lyapunov estimates}
    Let $X_0$ and $X_1$ be smooth vector fields on $\mathbb{R}^n$. Let $x_0$ be an isolated critical point of $X_0$ and $\mu >0$. Assume all eigenvalues of $JX_0(x_0)$ have real parts strictly smaller than $-\mu  $. Then there exists $\delta>0$ depending only on $X_0$ and a bounded neighborhood $U$ of $x_0$ as follows. Denoting by $\Phi^1_t$ the flow of $X_1$, we have that, if $\|X_0-X_1\|_{C^1}<\delta$, the the following hold:
    \begin{enumerate}[label=\textbf{\arabic*.}, leftmargin=*, align=left]
        \item[\textbf{1.} (Positive invariance)] For every $t\geq 0$, we have $\Phi^1_t(x)\in U$. In particular, $\Phi^1_t(x)$ is defined for every positive $t$;
        \item[\textbf{2.} (Uniqueness of critical point)] $X_1$ has a unique critical point $x_1$ in $U$; 
        \item[\textbf{3.} (Exponential estimates)] For every $m\geq 1$, there exists $c_m=c(m,X_0,\delta)>0$ such that, if $x\in U$ and $\|X_0-X_1\|_{C^1}<\delta$, then
    \begin{equation}\label{liapunov partial derivative}
         \|\Phi^1_t(x)-x_1\|\leq c_0e^{-\mu  t} \quad \text{and}\quad \|D^m\Phi^1_t(x)\|\leq c_me^{-\mu  t}(1+\|X_1\|_m),
    \end{equation}
    where $\|\cdot\|_m$ is the $C^m$ norm on $\overline{U}$. 
    \end{enumerate}
    A neighborhood $U$ with these properties is said to be $\Phi^1_t$-positively invariant. 
\end{lemma}
\begin{proof}
    See the appendix.
\end{proof}
We will also need the following technical lemma on matrices. Denote by $\sym^*_n$ the set of invertible symmetric matrices of size $n$.
\begin{lemma}\label{nonresonance lemma}
     Let $n\neq 2$ and $H$ be an invertible symmetric matrix of size $n$. Then there exists a matrix $F$ such that the eigenvalues of matrix
    \[
        B(H,F)=\frac{1}{2}\operatorname{tr}(F)H - FH,
    \]
    are all distinct and strictly negative. Moreover, there exists an open neighborhood $\mathcal H$ of $H$ in $\sym_n$ and a positive constant $\mu>0$ such that all the eigenvalues of $B(H',F)$ are strictly smaller than $-\mu$ for every $H'\in\mathcal H$.
\end{lemma}
\begin{proof}
    Let $H$ be a nonsingular symmetric matrix. We look for $F$ as in the statement that diagonalizes with $H$. In particular, we can assume that both $F$ and $H$ are diagonal. Let $F=\mathrm{diag}(f_1,\dots,f_n)$ and $H=\mathrm{diag}(h_1,\dots, h_n)$. In particular, the eigenvalues of $B(H,F)$ are given by
    \begin{equation}\label{sigma expression}
    \sigma_j=\frac12 h_j\sum_{i=1}^nf_i-h_jf_j.
    \end{equation}
    The vector $\sigma$ is obtained from $f$ by an invertible linear map. Indeed denote by $J$ be the $n\times n$ matrix where every entry is $1$. The matrix $\frac{1}{2}J-\id$ is invertible if and only if $n$ is different from $2$. Observing that $\sigma$ is given by
    \[
    \sigma=H\left(\frac12J-\mathrm{Id}\right)f
    \]
    proves the first part of the statement. The second part of the statement is given by the continuity of $B(\cdot,F)$ in the first entry.
\end{proof}
We now ready for the proof of  \cref{sol in crit}. We need to construct a solution of 
\[\dot\Delta_g[h](u)=v\]
in a neighborhood of the critical points of $u$. 
The strategy is the following. 
First, with \cref{jet solve}, we reduce to the case where $j^N_x v=0$ at each critical point $x$ of $u$. Then we reduce to a scalar equation using the ansatz $h=wF$, where $w$ is a  smooth compactly supported function and $F$ is a constant symmetric matrix.  The function $w$ will be found using the method of characteristics. Since the characteristics degenerate in a sink, we will need to prove that the constructed function is smooth. Tame estimates will follow naturally. 

\begin{proof}[proof of \cref{sol in crit} for $\dim(M)> 2$]
    As for the proof of \cref{jet solve}, we reduce the problem to a neighborhood of a single critical point $x_{u_0}$ of $u_0$. Let $(W,\varphi)$ be a geodesic normal chart with smooth boundary centered at $x_{u_0}$. For now, to simplify notation, we will treat $W$ as a bounded open subset of $\mathbb R^n$. Let $\mathcal V$ be a neighborhood of $q_0$. We can assume that $W$ contains a single critical point $x_u$ of $u$ for every $q=(g,u)$ in $\mathcal V$. We can further assume that every $x_u$ is contained in a fixed compact set $K\subseteq W$ with smooth boundary and nonempty interior. Let $F$ be a constant symmetric matrix. We look for a solution of \eqref{diff equation} in the form $wF$, where $w$ is a compactly supported smooth function. By \eqref{final 3.8}, we can express \eqref{diff equation} as
    \begin{equation}\label{simplified eq}
    a_qw+b_q(w)=v,
    \end{equation}
    where $a_q$ is a function, $b_q$ is a vector field and $v\in C^\infty(W)$. By \cref{nonresonance lemma}, we can pick $F$ such that $b_{q_0}$ is nondegenerate in $x_u$ and its Jacobian is negative definite. Up to shrinking $\mathcal V$, we can take a uniform constant $\mu>0$ such that the eigenvalues of $B_q=Db_q(x_u)$ are distinct and smaller than $-\mu$ for every $q\in\mathcal V$. Since we will mainly be interested in solving \eqref{simplified eq} in $K$, we can assume that $v\in C^{\infty}(\overline W)$ up to shrinking $W$. For now, we will simply treat these functions as functions on a subset of $\R^n$. Moreover, to relax notation, we will omit the dependence on $F$ unless relevant. Denote by $\Phi_{q,t}$ the flow of $b_q$. Up to a shrinking of $\mathcal V$, \cref{lyapunov estimates} gives a bounded neighborhood $U$ of $0$ which is $\Phi_{q,t}$--positively invariant for every $q\in \mathcal V$. We can assume $U=W$. Thus, for each integer $m\geq 1$, there exists a positive constant $c_m$ such that
    \begin{equation}\label{estimate to ref once}
        \|\Phi_{q,t}(x)-x_u\|\leq c_0e^{-\mu   t} \quad \text{and} \quad \|D^m\Phi_{q,t}(x)\|\le c_m e^{-\mu  t}(1+\|b_q\|_m).
    \end{equation}
    Moreover, $x_u$ is the only critical point of $b_q$ for every $q=(g,u)\in\mathcal V$. We can also assume that the $C^1$--norms of $a_q$ and $b_q$ are uniformly bounded by a positive constant $C$.
    \\
    We now solve equation \eqref{simplified eq} on $W$. First, assume that $v$ verifies $j^d_{x_u}(v)=0$ for a fixed integer $d>\frac{C}{\mu }$. Define
    $$I_q(t,x)\coloneqq\exp\left(\int_0^ta_q(\Phi_{q,\tau}(x))d\tau\right).$$
    We show that the following function is a well defined smooth solution of \eqref{simplified eq} in $W$ 
    \begin{equation}\label{solution}
        w_q(v)(x)=-\int_0^\infty v(\Phi_{q,t}(x))I_q(t,x) dt.
    \end{equation}
     
    \emph{Smoothness in $x$:} Let $m$ be a non-negative integer. Let $x\in W$. We want to show that
    \begin{equation}\label{integral to converge}
        \int_0^\infty \|D^{m}(v(\Phi_{q,t}(x))I_q(t,x))\| dt
    \end{equation}
    converges in $W$. This will prove that $w_q$ is a smooth function in $x$. From now on, we will denote by $C_m$ every constant depending only on $m,\mu  ,n,C$ and $W$. We will also denote the $C^j$ norms on $C^{\infty}(\overline{W})$ by simply $\|\cdot\|_{j}$. We have 
    \begin{equation}\label{main bound}
        \|D^{m}(I_q(t,x)v(\Phi_{q,t}(x)))\|\leq C_m\sum_{m_1+m_2=m}\|D^{m_1}I_q(t,x)\|\cdot\|D^{m_2}(v(\Phi_{q,t}(x))\|.
    \end{equation}
    We start by bounding the derivatives of $I_q$. First, observe 
    \[I_{q}(t,x)\leq e^{Ct}.\] We now assume $m_1\geq 1$. We have
    \begin{equation}\label{start inequality}
    \|D^{m_1}I_q(t,x)\|\leq C_m\sum_{k=1}^{m_1}I_q(t,x)\sum_{j_1+\ldots+ j_k=m_1} \int_0^{t}\|D^{j_1}G_{q,\tau}(x)\|d\tau\cdots\int_0^{t}\|D^{j_k}G_{q,\tau}(x)\|d\tau,
    \end{equation}
    where $G_{q,\tau}(x)=a_q(\Phi_{q,\tau}(x))$ and the $j$-indices are positive integers. Consider now the derivatives of $G_q$. Since composition of maps is a degree $1$ tame operation \cite[Lemma 2.3.4]{Hamilton}, we get
    $$
        \|D^{j}G_q(x,\tau)\|\leq C_m(1+\|\Phi_{q,\tau}\|_{j}+\|a_q\|_j)\leq C_m(1+e^{-\mu  \tau}\|b_q\|_j+\|a_q\|_j).
    $$
    Integrating, we obtain
    $$
        \int_0^t \|D^{j}G_q(x,\tau)\|d\tau\leq C_m\int_0^t (1+e^{-\mu  \tau}\|b_q\|_j+\|a_q\|_j)d\tau\leq C_m(t+\|b_q\|_j+t\|a_q\|_j).
    $$
    Thus,
    \begin{align}\label{use once}
        \int_0^{t}\|D^{j_1}G_{q,\tau}(x)\|d\tau\cdots\int_0^{t}\|D^{j_k}G_{q,\tau}(x)\|d\tau\leq C_m\prod_{i=1}^{k}(t+\|b_q\|_{j_i}+t\|a_q\|_{j_{i}})
    \end{align}
    We want to interpolate between the $C^j$ norms of $a_q$. We distinguish $2$ cases. First, assume $m_1=1$. Then $k=m_1=1$. Thus, there is a single $j$ to consider. Then \eqref{use once} becomes
    \[\int_0^{t}\|D^{1}G_{q,\tau}(x)\|d\tau\leq C_m(t+\|b_q\|_1+t\|a_q\|_{1})\leq C_mP_1(t)(1+\|b_q\|_1+\|a_q\|_1),
    \]
    where $P_1(t)=1+t$.\\
    \\
    Assume now $m_1>1$. Given a positive integer $i$ and a finite set $P$, we denote by $\binom Pk$ the collection of subsets of $i$ of cardinality $i$. If $P=\{1,\ldots, k\}$ is the collection of the first $k$ natural numbers, we'll write $\binom  Pk$ as $\binom ki$. We can write \eqref{use once} as
    \begin{align}\prod_{i=1}^k\int_0^{t}\|D^{j_i}G_{q,\tau}(x)\|d\tau&\leq C_m\prod_{i=1}^{k}(t+\|b_q\|_{j_i}+t\|a_q\|_{j_{i}})\nonumber\\&=C_m\left(t^k+\sum_{i=1}^{k}t^{k-i}\sum_{P\in\binom ki}\prod_{f\in P}(\|b_q\|_{j_f}+t\|a_q\|_{j_f})\right).
    \end{align}
    Fix $P\in\binom ki$. We have
    \begin{align}
        \prod_{f\in P}(\|b_q\|_{j_f}+t\|a_q\|_{j_f})&=\prod_{f\in P}\|b_q\|_{j_f}+\sum_{\ell=1}^it^\ell\sum_{P'\in\binom P\ell}\left(\prod_{\ell'\in P'}\|a_q\|_{j_{\ell'}}\cdot\prod_{\ell'\in\left(P'\right)^c}\|b_q\|_{j_\ell'}\right)
    \end{align}
    By \cref{interpolation}, we obtain
    \[\|a_q\|_{j}\leq C_m\|a_q\|_{1}^{\frac{m_1-j}{m_1-1}}\|a_q\|_{m_1}^{\frac{j-1}{m_1-1}}\leq C_m \|a_q\|_{m_1}^{\frac{j-1}{m_1-1}},\]
    where we have incorporated $\|a_q\|_1\leq C$ in the constant in the last step. The same inequality holds for the $j$--norms of $b_q$.
    This means that for $P'\in \binom P\ell$ we have
        \begin{equation}\label{estimate on aq}
            \prod_{\ell'\in P'}\|a_q\|_{j_{\ell'}}\prod_{\ell'\in\left(P'\right)^c}\|b_q\|_{j_\ell'}\leq C_m\|a_q\|_{m_1}^{\frac{-\ell+\sum_{\ell'\in P'} j_{\ell'}}{m_1-1}}\|b_q\|_{m_1}^{\frac{-i+\ell+\sum_{\ell'\in (P')^c} j_{\ell'}}{m_1-1}}=C_m \|a_q\|_{m_1}^{S(P')}\|b_q\|_{m_1}^{S((P')^c)},
        \end{equation}
        where $S(P')$ and $S((P')^c)$ are the corresponding exponents. We have the inequality
        \[S(P')+S((P')^c)\leq 1.\]
        Indeed, only set of indices $j$ with total sum smaller than $m_1$ appear in the summands for $S(P')$ and $S((P')^c)$. This means that
        \begin{equation*}
           (m_1-1)(S(P')+S((P')^c))=-\ell+\sum_{\ell'\in P'} j_{\ell'}-i+\ell+\sum_{\ell'\in (P')^c} j_{\ell'}=-i+\sum_{f\in P}j_f\leq m_1-1.
        \end{equation*}
        Then
        \[1-S(P')\geq S(P')+S((P')^c)-S(P')=S((P')^c).\]
        This gives
        \begin{align}
            \|a_q\|_{m_1}^{S(P')}\|b_q\|_{m_1}^{S((P')^c)}&\leq \|a_q\|_{m_1}^{S(P')}\left(1+\|b_q\|_{m_1}^{1-S(P'))}\right)\nonumber\\
            &\leq \|a_q\|_{m_1}^{S(P')}+C_m(\|a_q\|_{m_1}+\|b_q\|_{m_1})\nonumber\\
            &\leq 1+\|a_q\|_{m_1}+C_m(\|a_q\|_{m_1}+\|b_q\|_{m_1})\nonumber\\
            &\leq C_m(1+\|a_q\|_{m_1}+\|b_q\|_{m_1}),
        \end{align}
        where in the second inequality we have used
        \begin{equation*}
            \|a_q\|_{m_1}^{S(P')}\|b_q\|_{m_1}^{1-S(P')}\leq C_m(\|a_q\|_{m_1}+\|b_q\|_{m_1}).
        \end{equation*}
        This is a consequence of the convexity of $-\mathrm{log}$. Plugging everything into equation \eqref{start inequality}, we obtain
    \begin{equation}\label{first inequality}
        \|D^{m_1}I_q(t,x)\|\leq C_m(1+\|a_q\|_{m_1}+\|b_q\|_{m_1})P_{m_1}(t)e^{tC},
    \end{equation}
    where $P_{m_1}(t)$ is a polynomial. The same inequality also holds for $m_1=0$.\\
    \\
    We now study the terms of type $\|D^{m_2}(v(\Phi_{q,t}(x))\|$. Assume $m_2\geq 1$. The case $m_2=0$ will be a direct consequence of the following computations. We have
    \[
        \|D^{m_2}(v(\Phi_{q,t}(x))\|\leq C_m\sum_{k=1}^{m_2}\sum_{j_1+\ldots+j_k=m_2}\|D^kv(\Phi_{q,t}(x))\|\cdot\|D^{j_1}
        \Phi_{q,t}(x)\|\cdots\|D^{j_k}
        \Phi_{q,t}(x)\|,
    \]
    where the $j$ are positive integers. Applying \cref{lyapunov estimates} gives
    \begin{equation}\label{setup second ineq}
        \|D^{m_2}(v\circ\Phi_{q,t})(x)\|\leq C_m\sum_{k=1}^{m_2}\sum_{j_1+\ldots+j_k=m_2}\|D^kv(\Phi_{q,t}(x))\|e^{-k\mu t}\prod_{i=1}^k(1+\|b_q\|_{j_i}).
    \end{equation}
    Consider the products in the right hand side \eqref{setup second ineq}. It is a sum of terms in the form $\|b_q\|_{j_{i_1}}\cdots\|b_q\|_{j_{i_\ell}}$, where the $j$ are positive integers with $j_{i_1}+\cdots+j_{i_\ell}\leq m_2-(k-\ell)$. The $k-\ell$ term is given because all the $k-\ell$ indices that don't appear in the summand are positive. Assume that $m_2$ is different from $1$. The case $m_2=1$ give the same result with easier computations. As in \eqref{estimate on aq}, we have
    \begin{equation}\label{summing on b}
        \|b_q\|_{j_{i_1}}\cdots\|b_q\|_{j_{i_\ell}}\leq C_m\|b_q\|_{m_2}^{\frac{-\ell+\sum j_f}{m_2-1}}\leq C_m\left(1+\|b_q\|_{m_2}^{\frac{m_2-k}{m_2-1}}\right),
    \end{equation}
    where the last inequality follows by $j_{i_1}+\cdots+j_{i_\ell}\leq m_2-(k-\ell)$. We now give estimates to the $k$-th derivative of $v$. By \cref{interpolation}, we obtain
    \begin{equation*}
        \|D^kv(\Phi_{q,t}(x))\|\leq\|v\|_{k}\leq C_m\|v\|_{1}^{\frac{m_2-k}{m_2-1}}\|v\|_{m_2}^{\frac{k-1}{m_2-1}}.
    \end{equation*}
    We are interested only in tameness and convergence of \eqref{integral to converge}. Thus, we can assume $\|v\|_1$ to be bounded by a constant. This gives
    \begin{equation}\label{v estimates}
        \|D^kv(\Phi_{q,t}(x))\|\leq C_m\|v\|_{m_2}^{\frac{k-1}{m_2-1}}.
    \end{equation}
    In particular,
    \begin{align}\label{big k}
       \|v\|_{m_2}^{\frac{k-1}{m_2-1}}\left(1+\|b_q\|_{m_2}^{\frac{m_2-k}{m_2-1}}\right)\leq& \|v\|_{m_2}^{\frac{k-1}{m_2-1}}+\|v\|_{m_2}^{\frac{k-1}{m_2-1}}\|b_q\|_{m_2}^{\frac{m_2-k}{m_2-1}}\nonumber\\ \leq& 1+\|v\|_{m_2}+\|v\|_{m_2}^{\frac{k-1}{m_2-1}}\|b_q\|_{m_2}^{\frac{m_2-k}{m_2-1}}\nonumber\\
        \leq& C_m(1+\|v\|_{m_2}+\|b_q\|_{m_2}),
    \end{align}
    where in the last equality we used 
    \[
    \|v\|_{m_2}^{\frac{k-1}{m_2-1}}\|b_q\|_{m_2}^{\frac{m_2-k}{m_2-1}}\leq C_m(\|v\|_{m_2}+\|b_q\|_{m_2}).
    \]
    The estimate \eqref{big k} gives a bound for the summands in \eqref{setup second ineq}. However, we need a different estimate for $k\leq d$. The condition $j^d_{x_u}(v)=0$ gives
        $$\|D^k v(x)\|\le C_m\|x-x_u\|^{d-k}\|v\|_d.$$
        We can assume the $d$-norm of $v$ to be bounded by a constant $C_m$. This is because we are only interested in tameness and convergence of \eqref{integral to converge}. Therefore,
    \begin{align}\label{small k}
        \|D^kv(\Phi_{q,t}(x))\|e^{-k\mu t}\leq& C_m\|v\|_{d}\|\Phi_{q,t}(x)-x_u\|^{d-k}e^{-k\mu t}\nonumber
        \\ \leq& C_me^{-d\mu t},
    \end{align}
    where the last inequality follows from \eqref{liapunov partial derivative}. Observe that the right hand side does not depend on $k$. Then, applying \eqref{big k} and \eqref{small k} to \eqref{setup second ineq} gives
    \begin{align}\label{second inequality}
        \|D^{m_2}(v\circ\Phi_{q,t})(x)\|&\leq C_m\left(de^{-d\mu t}(1+\|b_q\|_{m_2})+\sum_{k=d+1}^{m_2}e^{-k\mu t}(1+\|b_q\|_{m_2}+\|v\|_{m_2})\right)\nonumber\\
        &\leq C_m(1+\|b_q\|_{m_2}+\|v\|_{m_2})\left(de^{-d\mu t}+\sum_{k=d+1}^{m_2}e^{-k\mu t}\right).
    \end{align}
    The bound in \eqref{second inequality} also holds for $m_2=0$ by \eqref{small k}. Combining \eqref{first inequality} and \eqref{second inequality} gives the following bound for $\|D^{m_1}I_q(t,x)D^{m_2}(v\circ\Phi_{q,t})(x)\|$
    \begin{equation}\label{super bound}
        C_m(1+\|a_q\|_{m_1}+\|b_q\|_{m_1})(1+\|b_q\|_{m_2}+\|v\|_{m_2})P_{m_1}(t)\left(de^{(C-d\mu) t}+\sum_{k=d+1}^{m_2}e^{(C-k\mu) t}\right).
    \end{equation}
    Since $P_{m_1}(t)$ is a polynomial and all the exponents appearing in \eqref{super bound} are strictly negative, \eqref{integral to converge} converges.  The right hand side of \eqref{super bound} is integrable on $[0,+\infty)$. Therefore,
    by dominated convergence $w_q(v)$  is a smooth function on $W$. Moreover,
    \begin{equation}\label{setup tame}
        \|w_q(v)\|_{m}\leq C_m(1+\|a_q\|_{m_1}+\|b_q\|_{m_1})\left(1+\|b_q\|_{m_2}+\|v\|_{m_2}\right).
    \end{equation}
    This is obtained by simply integrating \eqref{super bound} in $t$. We estimate $(1+\|a_q\|_{m_1})(1+\|b_q\|_{m_2})$ further. By \cref{lem:metric variation of laplacian}, $a_q$ and $b_q$ are tame in $g$. Since $m_1+m_2=m$, applying \cref{interpolation} gives
    \begin{align*}
        \|a_q\|_{m_1}\|b_q\|_{m_2}\leq C_m(\|a_q\|_m+\|b_q\|_m)
    \end{align*}
    and analogue estimates for the other product of norms in \eqref{setup tame}.
    Thus, we have tame estimates
    \begin{equation}\label{final estimates}
    \|w_q(v)\|_{m}\leq C_m(1+\|a_q\|_{m}+\|b_q\|_{m}+\|v\|_{m})\leq C_m(1+\|q\|_{m+r}+\|v\|_{m}),
    \end{equation}
    where $r$ is a nonnegative integer. In the second inequality we used that $a_q$ and $b_q$ are tame in $q$ by \cref{lem:metric variation of laplacian}.
    \emph{$w_q(v)$ is a solution:} We now show that $w_q(v)$ verifies \eqref{simplified eq}. By the method of characteristics, a smooth function $w$ is solution of equation \eqref{simplified eq} if and only if
    $$w(x)=I_q(s,x)w(\Phi_{q,s}(x))-\int_0^{s}I(t,x)v(\Phi_{q,t}(x))dt, \quad \text{for each $s\ge 0$.}$$
    Note that this formula is well defined, since we are working in a $\Phi_q$--positively invariant set. Let $x\in W$ and $s\geq 0$. We have 
    \begin{align*}
    -w_q(v)(\Phi_{q,s}(x))&=\int_0^{\infty}\mathrm{exp}\left(\int_{0}^ta_q(\Phi_{q,\tau}(\Phi_{q,s}(x)))d\tau\right) v(\Phi_{q,t}(\Phi_{q,s}(x)))dt
    \\&=\int_0^{\infty}\mathrm{exp}\left(\int_{0}^ta_q(\Phi_{q,\tau+s}(x))d\tau\right)v(\Phi_{q,t+s}(x))dt
    \\&=\int_0^{\infty}\mathrm{exp}\left(\int_{s}^{t+s}a_q(\Phi_{q,\tau}(x))d\tau\right)v(\Phi_{q,t+s}(x))dt\\
    &=\int_s^{\infty}\mathrm{exp}\left(\int_{s}^{t}a_q(\Phi_{q,\tau}(x))d\tau\right)v(\Phi_{q,t}(x))dt
    \\&=\int_s^{\infty}\frac{I_q(t,x)}{I_q(s,x)}v(\Phi_{q,t}(x))dt.\end{align*}
    Thus, $w_q(v)(\Phi_s(x))I_q(s,x)=-\int_s^{\infty}I_q(t,x)v(\Phi_{q,t}(x))dt$. Then $w_q(v)$ solves \eqref{simplified eq}.
    
\emph{Smoothness in the parameter:} It only remains to prove that $w_q(v)$ is smooth in $q$. We obtain this as an application of \cref{conti smooth}. Define $E_d\subseteq C^\infty(\overline W)$ as the set of smooth functions with $d$--th jet in the origin equal to $0$. Since $E_d$ is a closed vector subspace of $C^\infty(\overline W)$ of finite codimension, it is a tame Fr\'echet space. For every $q=(g,u)\in\mathcal V$ and smooth function $f$, we denote $T_{x_u}$ as $T_{x_u}f(x)=f(x+x_u)$. By \cref{lem:persistence of critical points}, this is a smooth tame operator. We denote by $\tilde a_q(x)=T_{x_u}a_q$ and by $\tilde b_q=T_{x_u}b_q$. Up to shrinking the neighborhoods we are working with, we can assume that they are well defined and all the estimates we have proven so far also hold the translated functions we'll use. For example, we could shrink $\mathcal V$ and take $W_1\subset W_2\subset W$ s.t. $W_1+x_u\subset W_2$, where $W_2$ is a $\phi_{q,t}$-positive invariant neighborhood of $W_1$ such that $W_2-x_u\subseteq W$ for every $(g,u)\in\mathcal V$. Consider the linear map
\[L:\mathcal V\times E_d\to E_d\]
defined as $L_q(v)=\tilde a_qv+\tilde b_q(v)$. Since $\tilde b_q(0)=b_q(x_u)=0$, this map is well defined i.e. $j^d_0L_q(v)=0$ for every $(q,v)$ in the domain of $L$. Recall from \cref{lem:persistence of critical points} and \cref{lem:metric variation of laplacian} that $x_u$, $a_q$ and $b_q$ are smooth tame in $q$. Thus, the map $L$ is also smooth and tame. Let now $q\in \mathcal V$. Consider the corresponding fiber
\[L_q:E_d\to E_d.\]
We define
\begin{equation}\label{last one before final check}
    \tilde w_q(v)(x)=-\int_0^\infty v(\Phi_{ q,t}(x+x_u)-x_u)I_{q}(t,x+x_u) dt
\end{equation}
As we argued in the first part of the proof, we have
\[D^m\tilde w_q(v)(x)=-\int_0^\infty D^m\left(v(\Phi_{ q,t}(x+x_u)-x_u)I_{q}(t,x+x_u)\right) dt\]
for every integer $m$. In particular, for $x=0$ the integrand is sum of terms containing a derivative of $v$ in $0$ of order less or equal than $m$. These derivatives are all $0$ for $m\leq d$, because $v$ belongs to $E_d$. Thus, $\tilde w_q(v)\in E_d$ for every $v\in E_d$. Moreover, $w_q(v)$ solves $\eqref{simplified eq}$, which implies that $\tilde w_q$ is a linear right inverse of $L_q$. Thus, we have shown that $L_q$ is surjective. We now show that it is injective. Since it is linear, it is enough to prove that $L_q(v)=0$ implies $v=0$. Since the $d$--th jet of $v$ in $0$ is null, there exists a positive constant $C_d$ such that
\[\|v(x)\|\leq C_d\|x\|^d\]
for every $x\in \overline{W}$. Moreover, by the method of the characteristics, $L_q(v)=0$ implies that
\[v(x)=v(\Phi_{q,s}(x+x_u)-x_u)I_{q}(s,x+x_u)\quad\text{for each }s\geq 0,x\in \overline W.\]
\begin{align*}
    |v(x)|&\leq |I_{q}(s,x+x_u)|\cdot|v(\Phi_{q,s}(x+x_u)-x_u)|\\
    &\leq e^{Cs}C_d\|\Phi_{q,s}(x+x_u)-x_u\|^d\\
    &\leq e^{Cs}C_de^{-d\mu s}=C_de^{(C-d\mu) s},
    \end{align*}
    which goes to $0$ as $s$ goes to infinity for $d>\frac C\mu$. Thus, $v$ is constantly $0$, implying injectivity of $L_q$. Moreover, the family $\tilde w_q$ of inverses tame by \eqref{final estimates}. It only remains to prove that it is also continuous. Let $(q_n,v_n)\in\mathcal V\times E_d$ be a sequence converging in $(q,v)\in\mathcal V\times E_d$. By the continuous dependence of the flow in the vector field, $\Phi_{q_n,t}(\cdot+x_{u_n})$ converges to $\Phi_{q,t}(\cdot+x_u)$ in the smooth topology. Then the corresponding integrand in \eqref{last one before final check} converges in the smooth topology uniformly on every bounded interval. Moreover, \eqref{super bound} applied to the translated functions gives an uniform integrable bound for the integrands of $\eqref{last one before final check}$ corresponding to the sequence $(q_n,v_n)$. Then, by dominated convergence we have $\tilde w_{q_n}(v_n)\to \tilde w_q(v)$.
    
     We can then apply \cref{conti smooth}. This gives that $\tilde w:\mathcal V\times E_d\to E_d$ is a smooth map. Observe now that $w_q(v)=T_{-x_u}\tilde w_q(T_{x_u}(v))$ for every $(q,v)$ such that $j^d_{x_u}v=0$. Thus, $w_q$ is smooth in $(q,v)$.
    
\emph{Definition of $R_{\crit}$:} We can finally defined $R_{\crit}$. Let $\{x_1,\ldots,x_m\}$ be the critical points of $u_0$. Let $d$ be sufficiently large. Consider sufficiently small $g_0$ geodesic charts $\varphi_i: W_i\to U_i\subset \R^n$ around them and for each one the corresponding compact set $K_i$ such that $K_i$ contains exactly one critical point of $u$ for every $(g,u)\in\mathcal V$. Let $w_{i,q}(v):U_i\to \R$ be the corresponding solutions of \eqref{simplified eq} obtained as above, each corresponding to a matrix $F_i$. Define the map
    \[R_{\crit}:\mathcal V\times C^{\infty}(M)\to \cG\]
    as
    \[R_{\crit}(q,v)=R_{\mathrm{jet},d}(q,v)+\sum_{i}\rho_i\varphi_{i,*}\left(w_{i,q}\left(\varphi_{i,*}^{-1}\left(v-\dot\Delta R_{\mathrm{jet},d}(q,v)\right)\right)F\right),\]
    where $R_{\mathrm{jet},d}$ is defined as in \cref{jet solve}. Since pullbacks are smooth tame operations by \cite[Lemma 2.3.4]{Hamilton}, this is a smooth tame map. The map $R_{\crit}$ is linear in $v$. Moreover, it verifies \eqref{new to cref} with $W_{\crit}=\cup_i \mathrm{int}(K_i)$ by construction. Indeed, for $x\in K_i$ we have
    \begin{align*}
        \dot\Delta R_{\crit}(q,v)(x)&=\dot\Delta R_{\mathrm{jet},d}(q,v)(x)+\dot\Delta \varphi_{i,*}\left(w_{i,q}\left(\varphi_{i,*}^{-1}\left(v-\dot\Delta R_{\mathrm{jet},d}(q,v)\right)\right)F_i\right)\\
        &=\dot\Delta R_{\mathrm{jet},d}(q,v)(x)+v-\dot\Delta R_{\mathrm{jet},d}(q,v)=v,
    \end{align*}
    where in the second equality we have used that $w_q$ solves \eqref{simplified eq} induced by the matrix $F_i$.
\end{proof}
A key step in the proof of \cref{sol in crit} was the fixing of a suitable matrix $F$ with \cref{nonresonance lemma}. This is not always possible in dimension $2$. However, the proof of \cref{sol in crit} is significantly simpler in dimension this case. Indeed, we can find an explicit solution using \cref{lem:metric variation of laplacian}.
\begin{proof}[proof of \cref{sol in crit} in $\dim(M)=2$]\label{dim 2}
    Let $\mathcal V$ be a neighborhood of $q_{0}$ and $W_{\crit}$ open in $M$ as follows. For every $(g,u)\in\mathcal V$, the set $\crit (u)$ is contained in $W_{\crit}$. Moreover, $\Delta_gu$ is never $0$ in $W$, where $W$ is a fixed open neighborhood of $W_{\crit}$. Let $\rho$ be a smooth function compactly supported in $W$ which is $1$ on $W_{\crit}$. Let $v$ be a smooth function on $M$. We look for a conformal solution of \eqref{diff equation} for $\rho v$. By \cref{lem:metric variation of laplacian}, this amounts to solving
    \[-\psi\Delta_g u=\rho v.\]
    We simply define $R_{\crit}(q,v)=-\frac{\rho v}{\Delta_g u}g$. This is a smooth tame map verifying the required property.
\end{proof}
\subsection{Right inverse of the linearized Laplace--Beltrami operator at regular points}
\cref{sol in crit} gives a smooth tame solution of \eqref{diff equation} in a fixed open neighborhood of critical points of $u$. We want to prove an analogue of \cref{sol in crit} outside such a neighborhood. Let $W\subseteq M$ be open connected with smooth boundary. We denote by $C^\infty_0(W)$ the space of smooth function vanishing at $\partial W$ with all derivatives. Similarly, we denote by $\mathfrak{X}_0(W)$ and by $\Gamma^{\infty}_0(S^2T^*W)$ the set of smooth vector fields and smooth $(0,2)$-tensors vanishing at $\partial W$ with all derivatives. On the set $\mathfrak X_0(W)$ we define $C^{k}$ norms as we did for smooth metrics in the introduction. This gives both vector spaces a structure of tame Fr\'echet spaces.
\begin{proposition}[Solution away from critical points]\label{sol away from crit}
    Let $q_0=(g_0,u_0)\in \mathcal N$. Then there exists a neighborhood $\mathcal{V}$ of $q_0$ in $\cG_+\times C^\infty(M)$ as follows. Let $W_{\mathrm{reg}}\subseteq M$ be an open connected set with smooth boundary, $\mathrm{crit}(u)\cap \overline{W}_{\mathrm{reg}}=\emptyset$ for every $(g,u)\in \mathcal V$, and $\Delta_{g_0}u_0$ is not identically $0$ in $\overline W_{\mathrm{reg}}$. Then there exists a smooth tame map
    \[R_\mathrm{reg}\colon \mathcal V\times C_0^\infty(W_{\mathrm{reg}})\to \Gamma^{\infty}_0(S^2T^*M)\]
    such that $\dot{\Delta}\circ R_{\mathrm{reg}}=\mathrm{Id}$.
\end{proposition}
It is very important to work with objects vanishing at the boundary. This is because $C_0^\infty(W_{\mathrm{reg}})$ can be embedded in $C^\infty(M)$, while $C^\infty(W_{\mathrm{reg}})$ cannot. The same holds for smooth metrics. This allows us to combine the solutions obtained from \cref{sol in crit} and \cref{sol away from crit} appropriately. By using \cref{lem:metric variation of laplacian}, we will reduce proving \cref{sol away from crit} to solving an equation of the type 
\begin{equation} \label{divergence}
\dive_g X=\psi
\end{equation}
where $\psi$ is a function supported on $W_{\mathrm{reg}}$ and $X$ is an unknown vector field. A solution of \eqref{divergence} supported again on $W_{\mathrm{reg}}$, is given by the Bogovskii operator. The following is a simplified version of \cite[Theorem III.3.3]{Galdi}.
\begin{theorem}[Bogovskii operator]\label{Bogovskii local}
    Let $U\subseteq \R^n$ be open, connected and bounded with Lipschitz boundary. Assume that $n\ge 2$. Then for every $f\in C^\infty_0(U)$ satisfying $\int_Uf=0$ there exists $\mathcal B(f)\in \mathfrak X_0(U)$ such that 
    $\dive \mathcal B(f)=f$. Moreover, $\mathcal B(f)$ can be chosen to be linear and such that
    \begin{equation}\label{Bogovskii tame estimate}
        \|\mathcal B(f)\|_{H_{l+1}(U)}\le C(\Omega,\ell)\|f\|_{H_l(U)}.
    \end{equation}
    for every positive integer $\ell$, where $C(\Omega,\ell)$ is a positive constants depending only on $\Omega$ and $\ell$.
    \end{theorem}
We need a generalization of this result for manifold with boundaries. The boundary assumption cannot be dropped in our method. This is because we aim to apply this result to $M\setminus W_{\crit}$, with $W_{\crit}$ coming from \cref{sol in crit}. This manifold always has a nonempty boundary. After making this generalization, we will discuss the tameness of the constructed Bogovskii operator. First, we define the appropriate spaces where this operator should be defined. Let $W$ be a compact connected smooth manifold with boundary. We denote the set of couples of metric and function with $0$ mean by
\[Z=\left\{(g,f)\in\cG_+\times C^{\infty}_0(W):\int_Wf\mathrm{d}\vol_g=0\right\}.\]
We will denote its fibers by
\[
Z(g)=\left\{f\in C^{\infty}_0(W)\mid(g,f)\in Z\right\}.
\]
We briefly argue for the tameness of these spaces.
\begin{lemma}[Tameness of $0$ mean]\label{0 mean}
    The sets $Z$ and $Z(g)$ are smooth tame Fr\'echet submanifolds.
\end{lemma}
\begin{proof}
    The set $Z(g)$ is a closed finite codimension vector subspace of $C^{\infty}_0(W)$. Thus, it is tame Fr\'echet. Fix a smooth metric $g_0$ on $W$. We denote by $\rho_g$ the unique positive smooth function such that $\mathrm{d}\vol_{g_0}=\rho_g\mathrm{d}\vol_g$. We first claim that the map
    \[\rho:\Gamma^\infty(S^2_+T^*W)\to C^{\infty}(W)\]
    defined as $\rho(g)=\rho_g$ is smooth tame. Indeed, by definition we have $\rho_g=\frac{\sqrt{\mathrm{det}(g_0)}}{\sqrt{\mathrm{det}(g)}}$. Thus, it is enough to prove that the map assigning a metric $g$ to the function $\mathrm{det}(g)$ is smooth tame. This immediately follows from the definition of $C^k$ norms. Consider now the map
    \[
    \chi :\cG_+\times C^{\infty}_0(W)\to \cG_+\times C^{\infty}_0(W)
    \]
    defined as $\chi(g,f)=(g,\rho_g f)$. This is a smooth tame embedding. By the definition of $\rho_g$, we have $Z=\chi(\cG_+\times Z(g_0))$.
\end{proof}
We have proven that the domains of the Bogovskii operators we will construct are tame. We will now define such operators. First consider the case where the metric is fixed.
\begin{theorem}[Global Bogovskii operator]\label{Bogovskii operator}
    Let $(W,g)$ be a compact connected Riemannian manifold with boundary. Assume $\dim W\ge 2$. Then, there exists a smooth tame linear map
    \[\mathcal{B}_g:Z(g)\to \mathfrak{X}_0(W)\]
    such that $\dive_g \mathcal{B}_g(f)=f$. In particular, for every positive integer $\ell$ there exists a positive constant $C(W,g,\ell)$ depending only on $W$, $g$ and $\ell$ such that
    \[\|\mathcal{B}_g(f)\|_{H_{\ell+1}(W)}\le C\|f\|_{H_\ell(W)}.\]
\end{theorem}
\begin{proof}
 Choose a finite coordinate cover
$\{(U_i,\varphi_i)\}_{i=1}^N$ such that each
$V_i=\varphi_i(U_i)$ is either a unit ball or a half unit ball. As in the proof of \cite[Lemma III.3.4]{Galdi}, we can find functions $f_i$ satisfying
    \begin{itemize}
        \item $f_i\in C^\infty_0(U_i)$,
        \item $\int_{U_i} f_i d \vol_g=0$,
        \item $\sum_i f_i=f$.
        \item $\|f_i\|_{H_l(W)}\le c_1\|f\|_{H_l(W)}$ for some positive constant $c_1=c_1(W,\ell)$.
    \end{itemize}
    Moreover, we can choose the $f_i$ to be linear in $f$. We solve the divergence equation in local coordinates. It becomes:
    \begin{equation}\label{divergence f_i on riemannian}
    \frac{1}{\sqrt{|g|}}\dive \left(\sqrt{|g|}v_i\right)=f_i,\end{equation} where $\mathrm{div}$ is the Euclidean divergence.
Multiplying both sides by $\sqrt{|g|}$ and then setting $w_i=\sqrt{|g|}v_i$ gives $\dive w_i= \sqrt{|g|}f_i$. The condition $\int_{U_i}f_id \vol_g=0$ in local coordinates becomes $\int_{V_i}f_i\sqrt{|g|}=0$. By \cref{Bogovskii local}, there exists $X_i(f_i)\in\mathfrak{X}_0(W)$ satisfying \eqref{divergence f_i on riemannian}. Since multiplication by $\sqrt{|g|}$ and its inverse are bounded on Sobolev spaces, we also have
\[\|X_i(f_i)\|_{H_{\ell+1}(W)}\le C_i\|f_i\|_{H_\ell(W)}.\]
Consider $\mathcal B_g(f)=\sum_iX_i(f_i)$. We have $\dive_g \mathcal B_g(f)=f$ by construction. Moreover,
\[\|\mathcal B_g(f)\|_{H_{\ell+1}(W)}\leq \sum_i\|X_i(f_i)\|_{H_{\ell+1}(W)}\leq C(W,g,\ell)\|f\|_{H_\ell(W)},\]
where $C$ is a positive constant depending only on $\ell, W$ and $g$.
\end{proof}

\begin{remark}[Regularity improvement]
    Observe that the Bogovskii operator is more than just tame. Indeed, it improves regularity. For our purposes, tameness will be enough.
\end{remark}
\cref{Bogovskii operator} is not sufficient for us. This is because we need to build a local right inverse of $\dot{\Delta}_g$. Thus, we must take into account a variation of the underlying metric. This is obtained in the following theorem.
\begin{theorem}[Augmented Bogovskii operator]\label{tameness of Bog}
    Let $(W,g_0)$ be a compact connected Riemannian manifold with boundary. Assume $\dim W\ge 2$. Then, the map
    \[
    \mathcal B(g_{0}):Z\to \mathfrak{X}_0(W)
    \]
    defined as $\mathcal B(g_0)(g,f)=\mathcal \rho_gB_{g_0}(\rho_g^{-1}f)$ is smooth and tame, where $\rho_g$ is the unique positive smooth function such that $\mathrm{d}\vol_{g_0}=\rho_g\mathrm{d}\vol_g$. Moreover, \begin{equation}\label{divergence condition}
    \dive_g\mathcal{B}(g_0)(g,f)=f.
    \end{equation}
\end{theorem}
\begin{proof}
    Observe that this map is well defined. Indeed,
    \[\int_Wf\mathrm{d}\vol_g=\int_W f\rho_g^{-1}\mathrm{d}\vol_{g_0}=0,\]
    where in the last line we used that $f$ belongs to $Z$.
    We will now prove tameness of $\mathcal B(g_0)$. Recall that the map $g\to \rho_g$ is smooth tame. This was proven in \cref{0 mean}. Thus, we can restrict proving tameness of $\mathcal B(g_0)$ to an open neighborhood of $(g_0,f_0)$ where the $1$-norm of $\rho_g$ and $\rho_g^{-1}$ are bounded by a fixed constant. This is because tameness is a local property. In this neighborhood, multiplication by $\rho_g$ and $\rho_g^{-1}$ are smooth tame maps by \cite[Lemma 2.3.4]{Hamilton}. Observe that applying a fixed linear continuous tame map such as $\mathcal B_{g_0}$ preserves smoothness and tame estimates up to a fixed constant. Thus, $\mathcal B(g_0)$ is smooth tame. Finally, we prove \eqref{divergence condition}. Observe that for every smooth vector field $X$ on $M$ and smooth metric $g$ the following holds \[
    \dive_g(X)=\rho_g\dive_{g_0}(\rho_g^{-1} X).
    \]
    This is obtained directly from the definition of divergence. Thus, for every $(g,f)\in Z$ we have
    \[\dive_g(\mathcal B(g_0)(g,f))=\rho_g\dive_{g_0}(\rho_g^{-1}\rho_g\mathcal B_{g_0}(\rho_g^{-1}f))=f.\]
\end{proof}
We are finally ready to prove \cref{sol away from crit}.
\begin{proof}[Proof of \cref{sol away from crit}]
    
    We divide the map $R_{\mathrm{reg}}$ as a sum of a conformal variation and a traceless variation, i.e.:
    \[R_{\mathrm{reg}}(q,v)= R_1(q,v)+ R_2(q,v),\]
     where $R_1$ is a conformal variation and $R_2 $ is a  traceless variation, both with respect to $g$.
     
    Let $y_0\in W_{\mathrm{reg}}$ be such $\Delta_{g_0}u_0(y_0)\neq 0$. Then there exists a neighborhood $V$ of $y$ in $W_{\mathrm{reg}}$ such that $\Delta_{g}u(y)\neq 0$ for every $(g,u)\in\mathcal V$ and $y\in V$. Let $\rho\geq 0$ be a nonzero bump function supported in $V$. Consider a conformal variation in the form $R_1(q,v)=c(v)\rho \Delta_g u\cdot g$, where $c=c(q,v)\in \R$. We look for $c$ such that
    \begin{equation}\label{norm variation} \int_{W_\mathrm{reg}} \dot{\Delta}_g[R_1](u) \, d\vol_g=\int_{W_\mathrm{reg}} v \, d \vol_g.\end{equation}
    By \cref{lem:metric variation of laplacian}, we get $\dot{\Delta}_g[R_1](u)=-c\rho (\Delta_gu)^2+\frac{n-2}{2}\langle \nabla (c\rho \Delta_gu),\nabla u\rangle_g$. Thus, integrating by parts gives
    \begin{align*}
        \int_{W_\mathrm{reg}} \dot{\Delta}_g[R_1](u) \, d\vol_g=&\int_{W_\mathrm{reg}} -c\rho (\Delta_gu)^2+\frac{n-2}{2}\langle \nabla (c\rho \Delta_gu),\nabla u\rangle_g\, d\vol_g \\
        =&\int_{W_\mathrm{reg}} -c\rho (\Delta_gu)^2-\frac{n-2}{2} c\rho (\Delta_gu)^2\, d\vol_g\\
        =& -\frac n2 c\int_{W_\mathrm{reg}} \rho (\Delta_gu)^2 \, d\vol_g
    \end{align*}
    In particular, $ \rho (\Delta_gu)^2$ is a nonnegative function which is positive on an open subset of $W_{\mathrm{reg}}$. Therefore, we can choose 
    \[c=-\frac{2\int_{W_{\mathrm{reg}}}vd\vol_g}{n\int_{W_\mathrm{reg}} \rho (\Delta_gu)^2 \, d\vol_g}\]
    so that \eqref{norm variation} holds. It follows also that the map $R_1$ defined in this way is smooth and tame, being composition of integrals and multiplication by smooth functions.

Let $\tilde v\coloneqq v-\dot \Delta_g[ R_1((q,v)](u)$. We now search for a traceless variation $R_2(q,\tilde v)$ such that 
\[\dot{\Delta}_g[R_2(q,\tilde v)](u)=\tilde v.\]
By \cref{lem:metric variation of laplacian}, this amounts to solve
\begin{equation}\label{divergence equation}
-\dive_g(Du \lrcorner R_2^*)=\tilde v,
\end{equation}
where $R^*_{2}$ is the tensor $R_{2}$ but all the indexes are raised by $g$.
By \cref{tameness of Bog}, we can write \eqref{divergence equation} as
\begin{equation}\label{traceless dive}Du\lrcorner R^*_{2}=\mathcal{B}(g_0)(g,-\tilde v).\end{equation}
Denote by $B\coloneqq \mathcal{B}(g_0)(g,-\tilde v)\in \mathfrak{X}_0(W_{\mathrm{reg}})$ and by $B^\perp\coloneqq B-\frac{\langle \nabla u,B\rangle_g}{\|\nabla u\|^2}\nabla u$ its projection on $(\nabla u)^\perp$. Notice that $B^\perp$ is a smooth vector field because $\nabla u$ never vanishes on $W_{\mathrm{reg}}$. Define $R^*_{2}$ as:
\begin{equation}\label{h_2 expression}R^*_{2}\coloneqq \frac{\nabla u \otimes B^\perp+B^\perp\otimes\nabla u}{\|\nabla u\|^2}+\frac n{n-1}\frac{\langle B,\nabla u\rangle_g}{\|\nabla u\|^2}\left(\frac{\nabla u \otimes \nabla u}{\|\nabla u\|^2}-\frac{1}{n}g^*\right)\end{equation}
    Since $V\equiv 0$ in a neighborhood of the critical points of $u$ then $R^*_{2}$ is smooth.
    Notice that $R^*_{2}$ is symmetric. A simple computation shows that it is traceless:
    \[\tr R^*_{2}=n\frac{\langle B,\nabla u\rangle}{\|\nabla u\|^2}\left(1-1\right)=0.\]
    Also \eqref{traceless dive} is satisfied:
    \begin{align*}Du\lrcorner R^*_{2}&= B^\perp -\cancel{Du(B^\perp)\frac{\nabla u}{\|\nabla u\|^2}}-n\frac{\langle B,\nabla u\rangle_g}{\|\nabla u\|^2}\left(-\nabla u+\frac{1}{n}\nabla u\right)\\
    &=B^\perp +\frac{\langle B,\nabla u\rangle_g}{\|\nabla u\|^2} \nabla u\\
    &= B. \end{align*}
Since the Bogovskii operator is smooth and tame by \cref{tameness of Bog}, it follows that
$$R_2\colon \mathcal V \times C^\infty_0(W_{\mathrm{reg}})\to \Gamma^\infty_0(S^2T^*W_{\mathrm{reg}})\subseteq \Gamma^\infty(S^{2}T^*M)$$ is a smooth tame map since it is a composition of smooth tame maps. Finally we define $R_{\mathrm{reg}}=R_1+R_2$, so that 
\begin{align*}
    \dot \Delta_g[R_{\mathrm{reg}}((g,u),v)](u)&= \dot \Delta_g[R_1((g,u),v)](u)+\dot \Delta_g[R_2(g,u),v)](u)\\
    &=\dot \Delta_g[R_1((g,u),v)](u)+v-\dot \Delta_g[ R_1((g,u),v)](u)
    \\ &= v.
\end{align*}
\end{proof}
\begin{remark}
    Contrary to the one of $R_\mathrm{{crit}}$, the smoothness in $q$ of $R_{\mathrm{reg}}$ can be easily derived from its explicit expression.
\end{remark}
\subsection{Density of the submersivity set}
So far, we have essentially constructed the map $R$ in the statement of \cref{laplacian use}. To prove the theorem, it is still necessary to show that $\mathcal N$ is open dense. We will need a technical transversality result to prove do so. The following can be found in \cite[17.2, 19.1, 19.2]{transversalmapping}.

\begin{theorem}\label{transversality theorem}
Let $X$ and $Y$ be $C^1$ Banach manifolds with a countable base. let $W\subset Y$ be $C^1$ Banach submanifold of finite codimension $c$, and let $M$ be a smooth manifold of dimension $n$. Consider a $C^1$map
$
F \colon X \times M \to Y,
$
and assume that $n<c<+\infty$. If $F\pitchfork W$,  then denoting $F_x\coloneqq F(x,\cdot):M\to Y$, we have that the set
\[
\bigg\{x\in X \;\bigg|\; \operatorname{Im}(F_x)\cap W=\emptyset \bigg\}
\]
is residual in $X$.
\end{theorem}
\begin{proof}
    The proof is a consequence of standard transversality theory in Banach manifolds. Since this result is standard and is not used elsewhere in the present work, we do not recall the definitions of transversality and Fredholm maps for Banach manifolds. We only give a sketch of the proof, and we refer the interested reader to \cite{transversalmapping} for the precise theory and definitions.
    
    Since $F\pitchfork W$, from \cite[Corollary 17.2]{transversalmapping} we obtain that $F^{-1}(W)$ is a $C^1$ Banach submanifold of $X\times M$ of codimension $c$. Consider the restriction of the projection on the $X$ component
    \[\pi|_{F^{-1}(W)}\colon F^{-1}(W)\to X.\]
    Then by \cite[Lemma 19.2] {transversalmapping} $\pi|_{F^{-1}(W)}$ is a $C^1$ Fredholm map of index $n-c<0$. By the Sard--Smale theorem \cite[Theorem 16.2]{transversalmapping} the set of regular values of $\pi|_{F^{-1}(W)}$ is residual in $X$. We now observe that, since the Fredholm index is negative, a regular value of $\pi|_{F^{-1}(W)}$ cannot belong to its image. Thus the set of regular values is precisely \[ X\setminus \pi\bigl(F^{-1}(W)\bigr). \] 
    Finally,  $x\notin \pi\bigl(F^{-1}(W)\bigr)$  is equivalent to $F_x(M)\cap W=\emptyset.$ Hence \[ \bigg\{x\in X \;\bigg|\; \operatorname{Im}(F_x)\cap W=\emptyset \bigg\} \] is residual in $X$.
\end{proof}
To apply \cref{transversality theorem}, it will be convenient to work with $C^k$ spaces. The following topological lemma guarantees that this is possible.
\begin{lemma}\label{lemma Q_k residual}
Let $\{X_k\}_{k\ge k_0}$ be a nested family of Banach manifolds with continuous dense inclusions
\[
X_\ell\hookrightarrow X_k,\qquad \ell\ge k,
\]
and let
\[
X_\infty:=\bigcap_{k\ge k_0}X_k
\]
be endowed with the induced Fréchet topology. Assume that the inclusion $X_\infty\hookrightarrow X_k$ is dense for every $k$. Let $D\subset X_k$ and assume that
\[
D_\ell=D\cap X_\ell,\qquad \ell\ge k
\]
is open dense in $X_\ell$ for every sufficiently large $\ell$. Then
\[
D_\infty:=D\cap X_\infty
\]
is open dense in $X_\infty$.
\end{lemma}
 \begin{proof}
Openness is clear, since $X$ is endowed with the Fre\'chet topology induced by the limit of the $X_\ell$. It remains to prove density. Let $U\subset X_\infty$ be a nonempty open set. By definition of the Fr\'echet topology on $X_\infty$, there exists $\ell$ and a nonempty open set $V\subset X_\ell$ such that $V\cap X_\infty\subseteq U.$  We can assume that $\ell$ is large enough that $D_\ell$ is dense in $X_\ell$. Then, the set $V\cap D_\ell$ is a nonempty open subset of $X_\ell$. Since $X_\infty$ is dense in $X_\ell$, the intersection $X_\infty\cap V\cap D_\ell=V\cap D_\infty$ is nonempty. Since $V\cap X_\infty\subseteq U$, this completes the proof.
\end{proof} 
\begin{remark}[Generic properties and regularity]
    \cref{lemma Q_k residual} allows to work in the Banach setting when studying some open dense properties of smooth objects. However, not all open dense or residual properties of smooth functions come from finite regularity ones. This justifies working in the Fr\'echet setting.
\end{remark}
We can now prove that $\mathcal N$ is open dense.
\begin{proposition}\label{N open dense}
    The set $\mathcal N$ is open dense. In particular, it is a tame Fr\'echet variety.
\end{proposition}
\begin{proof}
     For $\dim(M)>2$, this result is equivalent to Morse functions being open dense in $C^\infty(M)$. Then we can assume $\dim(M)=2$. By \cref{lemma Q_k residual}, it is enough to prove that the set $\mathcal N_{C^k}$ is open dense for every $k$ sufficiently large. Since the set of Morse functions is open dense, it is enough to study the second condition. We first prove density as an application of \cref{transversality theorem}. Consider the map
    \begin{align*}
         F\colon 
         \left(\cG_+^{k}\times H_{k}(M)\right)\times M&\to \left(\cG_+^k\times H_{k}(M)\right) \times T^*M \oplus \R\\
         (g,u),x&\mapsto (g,u),  D_x u,\Delta_gu.
     \end{align*}
     Let $Z$ be the $0$ section of $ \left(\cG_+^{k}\times H_k(M)\right) \times T^*M \oplus \R$ and define $W=\left(\cG_+^{k}\times H_k(M)\right) \times Z$. Then $W$ is exactly the set we need to avoid generically. Since $\codim(W)=3> 2=\dim(M)$, we need only to prove that $F\pitchfork W$ to apply \cref{transversality theorem}. To show this, simply observe that we can always pick a smooth function on $M$ with a prescribed jet at a fixed point. Finally, $\mathcal N$ is open because $M$ is compact. Thus, the set of $q\in\cG_+^{k}\times H_k(M)$ such that $F(q,x)$ avoids $Z$ is open.
\end{proof}
We can finally prove \cref{laplacian use}. It is simply a combination of \cref{sol in crit} and \cref{sol away from crit}.
\begin{proof}[{Proof of \cref{laplacian use}}]
    We apply \cref{sol in crit} to $q_0$. This gives a neighborhood $\mathcal V_1$ of $q_0$, an open set $W_\mathrm{crit}\subseteq M$ and a map $R_\mathrm{crit}$ as in \cref{sol in crit}. Let $K\subseteq W_{\crit}$ be a compact set with smooth boundary containing $\mathrm{crit}(u_0)$ in its interior. Up to shrinking $\mathcal{V}_1$, we can assume that $\mathrm{crit}(u)$ is contained in $K$ for every $(g,u)$ in $\mathcal{V}_1$. Moreover, we can assume that $M\setminus K$ is connected since $n\geq 2$. This is because $\crit(u_0)$ is finite. Since $u_0$ is not constant, $\Delta_{g_0}u_0$ is nonzero in an open subset of $M$. Up to shrinking $\mathcal{V}_1$ and $K$, we can assume that it is not constantly zero in the complement of $K$.
    
    We can now apply \cref{sol away from crit} with $W_{\mathrm{reg}}=M\setminus K$. This gives a neighborhood $\mathcal V_2$ of $q_0$ and a map $R_{\mathrm{reg}}$. We can assume $\mathcal V\coloneqq\mathcal V_1=\mathcal V_2$. Consider the open set $W_\mathrm{reg}=M\setminus K$. We define
    \[R\colon\mathcal V\times C^\infty(M)\to \cG\]
    as $R(q,v)=R_\mathrm{crit}(q,v)+R_\mathrm{reg}(q,v-\dot{\Delta}R_\mathrm{crit}(v))$. This map is well defined. This is because $v-\dot{\Delta}R_\mathrm{crit}(v)\in C^\infty_0(W_\mathrm{reg})$ by \cref{sol in crit}. Moreover, it is a sum of smooth tame maps. Thus, it is smooth tame. Finally,
    \begin{align*}
        \dot{\Delta}R(q,v)=\dot{\Delta}R_\mathrm{crit}(q,v)+\dot{\Delta}R_{\mathrm{reg}}(v-\dot{\Delta}R_\mathrm{crit}(q,v))\\
        =\dot{\Delta}R_\mathrm{crit}(q,v)+v-\dot{\Delta}R_\mathrm{crit}(q,v)=v.
    \end{align*}
\end{proof}

\section{Generic properties of eigenfunctions of the Laplace--Beltrami operator}
\cref{laplacian use} can be interpreted as a sort of flexibility property for the Laplace--Beltrami operator. It states that every tangent direction can be realized by $\Delta$ only using a variation of the metric. In this section, we discuss the consequence of this result on the generic properties of the eigenfunctions of second order selfadjoint elliptic operators. We will call a set containing a countable intersection of open dense sets residual. We will say that a property is generic if it is verified by every element in a residual set. To even begin to discuss whether a property of selfadjoint elliptic operators is generic, it is necessary to equip them with a topology. In this case, the topology will be induced by a parametrization of the kind of elliptic operators we want to consider. First, we'll consider properties of operators of the Laplace--Beltrami type. They will be parametrized naturally with the identification $g\to\Delta_g$. This case represents the most direct consequence of \cref{laplacian use}. Afterwards, we'll consider operators of the same type, but parametrize them with metrics induced by embeddings of $M$ into the Euclidean space. Finally we'll consider consider the family of Schr\"odinger operators with independent metric and potential parameters.
\subsection{Genericity properties after Uhlenbeck}
Generic properties of eigenfunctions of the Laplace--Beltrami operator were previously studied by K. Uhlenbeck in \cite{Uhlenbeck} in the Banach setting. Her argument is carried out in a finite-regularity Banach setting and is based on the Sard--Smale theorem and parametric transversality for Fredholm maps between Banach manifolds. We briefly recall the part of her work that is relevant to our setting.

Denote by $\cG_+^k$ the space of $C^k$ sections on $C^k$ symmetric $(0,2)$--tensors on $M$. Each element of $\cG_+$ induces a Laplace--Beltrami operator $\Delta_g:H_{k}(M)\to H_{k-2}(M)$. Consider the following properties of element in $\cG_+^k$
\begin{enumerate}
    \item[(1)] Every eigenvalue of $\Delta_g$ is simple.
    \item[(2)] Every eigenfunction of $\Delta_g$ has $0$ as regular value.
    \item[(3)] Every eigenfunction of $\Delta_g$ is Morse.
\end{enumerate}
One can show the following.
\begin{theorem}[{\cite[Theorem 8]{Uhlenbeck}}]\label{Uhlenbeck}
    Let $M$ be a compact connected smooth manifold with $\dim(M)=n\geq 2$. Let $k>n+3$ be a positive integer. Then, the set 
    \[\{g\in \cG_+^k: \Delta_g\mathrm{\ satisfies\ (1),\ (2),\ (3)\ on\ nonconstant\ eigenfunctions}\}\]
    is residual in $\cG_+^k$.
\end{theorem}

\cref{Uhlenbeck} was originally proven using transversality techniques on Banach manifolds, which are not available in this form in the Fr\'echet setting. This justifies the different approach of our work, which instead relies on \cref{nash moser}. However, notice that \cref{Uhlenbeck} will play a key role in the proof of \cref{tame}. This density input is essential for the main theorem; its fixed-metric potential version is also used for the Schr\"odinger extension. In this sense Uhlenbeck's study of the Banach setting provide the groundwork to treat the smooth one.

\begin{remark}
    Observe that there is no analogue of \cref{Uhlenbeck} in dimension $1$. Indeed, every $1$-dimensional compact connected smooth Riemannian manifold is isometric to $S^1$ with a uniform metric and every nonzero eigenvalue of the corresponding Laplace--Beltrami operator has multiplicity $2$.
\end{remark}

Uhlenbeck's method also applies to other families of operator. We will use the following to study a special class of operators, called Schr\"odinger operators, in \cref{sect Schro}.
\begin{theorem}[{\cite[Theorem 7]{Uhlenbeck}}]\label{fiber uhlenbeck}
    Let $M$ be smooth compact connected manifold of any dimension and $k\geq 3$ be a positive integer. Denote by $L$ a fixed second order self--adjoint elliptic operator on $M$ with smooth coefficients. Then the set
    \[\{v\in C^k(M): L+v\mathrm{\ satisfies\ (1),\ (2),\ (3)}\}\]
    is residual in $C^k(M)$.
\end{theorem}
\begin{remark}[Eigenpair space formulation]
    Despite being different from our methods, the proofs of \cref{Uhlenbeck} and \cref{fiber uhlenbeck} still pass through the eigenpair space, although Sobolev spaces are used in the original work. In particular, Uhlenbeck's results could also be stated as density of properties (1) to (3) for triples in the eigenpair space. This formulation will be more convenient for us and will be implicitly used in many occasions.
\end{remark}
To conclude this section, we explain why only operators of order $2$ are considered in \cite{Uhlenbeck}. This is due to the the use the \emph{weak unique continuation theorem}, which is a key tool in the proofs of \cref{Uhlenbeck}.
\begin{theorem}[{\cite[Weak Unique Continuation Theorem]{Uhlenbeck}}]\label{weak}
    Let $L$ be a second order elliptic operator with $C^3$ coefficients. Let $D$ be a nonempty open connected subset of $M$. If $Lv=0$ in $D$ and $v=0$ on a nonempty open subset of $D$, then $v$ is identically $0$ in $D$.
\end{theorem}
Notice how \cref{weak} requires the order of the operator to be $2$. Indeed, ellipticity alone does not guarantee a unique continuation of solution and \cref{weak} it does not hold for higher order operators. Since the present paper requires \cref{Uhlenbeck}, we also restrict ourselves to studying operators of second order.

\subsection{Geometry of the eigenpair space}
Consider the eigenpair set
\[
\mathcal Q'=\bigg\{(g,\lambda,u)\in 
\cG_+\times \R_{>0} \times C^\infty(M)
\;\bigg|\;
(\Delta_g+\lambda)u=0,\ m_\lambda(g)=1, u\neq 0
\bigg\}.
\]
Near a simple eigenpair, the eigenvalue and the corresponding eigenfunction
are determined by the metric up to multiplication by a nonzero scalar. Thus one expects $\mathcal Q'$ to be locally parametrized by the metric together with one real normalization parameter. Our goal is to make this precise by showing that the space $\mathcal Q'$ is a tame Fréchet submanifold of the space
$\cG_+\times \mathbb R_{>0}\times C^\infty(M),$
modeled on $\cG\times \mathbb R$. The basic idea is to encode the eigenvalue equation and the normalization of the eigenfunction in a single map. Consider
\begin{align*}
\Psi\colon
\cG_+\times\mathbb R_{>0}\times C^\infty(M)&
\longrightarrow
\cG_+\times\mathbb R\times C^\infty(M),\\
(g,\lambda ,u)&\longmapsto\left( g,\frac 12\|u\|_{L^2_g}^2,(\Delta_g+\lambda)u\right)
\end{align*}
Let $(g,\lambda,u)$ be a point in $\mathcal Q'$ and assume that $-\lambda$ is the $k$--th eigenvalue of $\Delta_g$. Locally, the space $\mathcal Q'$ can be described as the inverse image under $\Psi$ of $U\times \mathbb R\times\{0\}$, where $U$ is an open neighborhood of $g$ in $\cG_+$. Thus, if $\Psi$ is locally invertible near a simple eigenpair, the desired manifold structure on $\mathcal Q'$ follows from the corresponding product structure in the target.

This reduces the problem to an application of the  inverse function theorem. Since the spaces
involved are Fréchet spaces, the ordinary Banach-space inverse function theorem does not apply. We will instead use the Nash--Moser inverse function theorem, i.e.
\cref{nash moser}. 

To prove the required tame estimates for the inverse of $D\Psi$, it is
convenient to use Sobolev norms. The analytic input is the following refinement of the elliptic estimates, with a bound also on the metric parameter.
 \begin{lemma}[Local elliptic estimates]\label{elliptic estimates}
     Let $q_0=(g_0,\lambda_0,u_0)\in\mathcal Q'$. Then there exists a neighborhood $\cV$ of $q_0$ and positive constants $\{C_n\}_{n\in\mathbb{N}}$ such that for all $q=(g,\lambda,u)\in \mathcal V$, $v\in C^\infty(M)$ and  $\eta \in \R$ we have
     \begin{align*}
        \|v\|_{H_2}+|\eta|\leq C_0\left(\|(\Delta_g+\lambda)v+\eta u_0\|_{H_0}+\left|\langle u_0,v\rangle_{L^2_{g_0}}\right|\right).
     \end{align*}
     Moreover, $\|v\|_{H_{n+2}}$ is bounded by
     \begin{align*}
     C_n\left(\|(\Delta_g+\lambda)v+\eta u_0\|_{H_{n}}+\|(\Delta_g+\lambda)v+\eta u_0)\|_{H_0}\|q\|_{n+1}+(1+\|q\|_{n+1})\left|\langle u_0,v\rangle_{L^2_{g_0}}\right|\right)
     \end{align*}
     for every $q=(g,\lambda,u)\in\cV$ and $n\in\mathbb{N}$.
\end{lemma}
\begin{proof}
      The following result is a particular case \cite[Lemma 3.3.4 Sec 2]{Hamilton}. With the notation of the original statement, it is obtained by considering the map $i(\dot\lambda)=\dot\lambda u_0$ and the map $j(v)=\langle u_0,v\rangle_{L^2_{g_0}}$.
\end{proof}
We are ready to show that $\mathcal Q'$ is a Fréchet submanifold.
\begin{theorem}[$\mathcal Q'$ is tame Fr\'echet]\label{simple}
 Let $q_0=(g_0,\lambda_0,u_0)\in\mathcal Q'$. Then, there exists a neighborhood $\mathcal V$ of $q_0$ such that $\mathcal Q'\cap\mathcal V$ is a smooth tame Fr\'echet submanifold of $\cG_+\times \R_{>0}\times C^\infty(M)$.
\end{theorem}
\begin{proof}
     Consider the smooth tame map between Fre\'chet manifolds
     \begin{equation}\label{Psi}
 \Psi\colon\cG_+\times\R_{>0}\times C^{\infty}(M)\to \cG_+\times\R\times  C^{\infty}(M),
\end{equation}
defined as $\Psi(g,\lambda,u)=\bigl(g,\ \frac12 \|u\|^2_{L^2_g},(\Delta_g+\lambda)u\bigr)$. Its differential in a point $q=(g,\lambda,u)$ can be written as
\[
D_{q}\Psi(h,\eta,v)=
 \Bigl(h,\ \frac14\int_M u^2\,\mathrm{tr}_{g}(h)\,d\mathrm{vol}_{g}+\langle u,v\rangle_{L^2_g},\dot \Delta_{g}(u)[h]+(\Delta_g+\lambda)v+\eta u\Bigr),
\]
where $\mathrm{tr}_{g}(h)$ is the trace with respect to $g$ and we used $d(\mathrm{vol}_{g})[h]=\tfrac12\mathrm{tr}_g(h)\,d\mathrm{vol}_g$. To apply \cref{nash moser}, we need to prove that the differential is locally invertible with a family of continuous inverses. This amounts to showing that the equation
\begin{equation}\label{correction equation}
    D_{q}\Psi(h,\eta,v)=(h',\eta',v')
\end{equation}
has a unique solution for every $h',\eta',v'$. Then, we need to prove that this family of solutions is smooth and tame in $q$. Let $\mathcal V$ be an open neighborhood of $q_0=(g_0,\lambda_0,u_0)$. If $\lambda_0$ is a simple $k$-th eigenvalue we can assume that all the metrics that appear in $\mathcal V$ have simple $k$--th eigenvalue by \cite[Theorem 2.2]{Bando}. We now argue that it is enough to study only the second and third component of $D\Psi$. Indeed, \eqref{correction equation} forces $h=h'$. Thus, the equation reduces to 
\begin{equation}\label{second correction equation}
    D_{q}\Psi(0,\eta,v)=D_q \Psi(h,\eta,v)-D_{q}\Psi(h,0,0)=(0,\eta'',v'').
\end{equation}
Observe that the first component of $D_{q}\Psi(h,\eta,v)$ depends only on $h$. In particular, requiring for it to be $0$ imposes no condition on $\eta$ and $v$. This implies also that the tameness of the solution map of \eqref{correction equation} follows from the  tameness of the solution map of \eqref{second correction equation}. Therefore we can reduce the problem to studying 
 \[
 D\tilde{\Psi}\colon\mathcal V\times \R\times C^{\infty}(M)\to \R\times C^{\infty}(M),\quad D_q\tilde{\Psi}(\eta, v)=\left(\langle u,v\rangle_{L^2_g},(\Delta_g+\lambda)v+\eta u\right).
 \]
 A simple computation shows that $D_{q_0}\tilde{\Psi}$ is invertible. We first prove that we can take $\mathcal V$ such that $D\tilde{\Psi}$ is invertible with a family of continuous inverses. Fix a sufficiently large positive integer $\ell$. Consider the extension
\[
D\tilde\Psi_{\ell}:\mathcal V\times\R\times H_{\ell+2}(M)\to \R\times H_{\ell}(M),
\]
defined as in the smooth case. The assignment $q$ to $D_q\tilde\Psi_{\ell}$ is a continuous map from $\mathcal V$ to $\mathcal L=\mathcal L(\R\times H_{\ell+2},\R\times  H_{\ell})$. Observe that $\mathcal L$ is the set of continuous maps between $2$ Banach spaces. In particular, the subset of invertible operators is open. Thus, we can take $\mathcal V$ such that $D_q\tilde\Psi_{\ell}$ is invertible for every $q$ in $\mathcal V$. Moreover, the total inverse 
\[D\tilde\Psi_{\ell}^{-1}\colon \mathcal V\times \R \times H_{\ell}(M)\to \R\times H_{\ell+2}(M)\]
is continuous. This is because taking inverses is a continuous operation in Banach spaces. By \cref{schauder estimates}, the restriction
\[D\tilde\Psi^{-1}\coloneqq D\tilde\Psi_{\ell}^{-1}|_{C^{\infty}(M)}\colon \mathcal V\times \R \times C^{\infty}(M)\to \R\times C^{\infty}(M)\]
is s.t. $D_q\tilde\Psi^{-1}$ is invertible for every $q\in \mathcal V$. We want to show that the total map $D\tilde\Psi^{-1}$ is continuous. Let $(q_n,c_n,f_n)$ be a sequence in $\mathcal V\times \R \times C^{\infty}(M)$ converging in the smooth topology to $(q,c,f)\in \mathcal V\times \R \times C^{\infty}(M)$. Denote by $D\tilde\Psi_{q_{n}}^{-1}(c_n,f_n)=(\eta_n,v_n)$. Continuity in the $H_\ell$ topology already gives $\eta_n$ converges to $\eta$ and $v_n$ to $v$, the latter only in the $H_\ell$ topology. Consider the PDEs
\[(\Delta_{g_n}+\lambda_n)v_n+\eta_n u_n=f_n\]
and
\[(\Delta_g+\lambda)v+\eta u=f.\]
By subtracting them, we obtain
\begin{equation}\label{one more}
    (\Delta_{g_n}+\lambda_n)(v_n-v)=(f_n-f)+(\eta_n u_n-\eta u)+\left((\Delta_g+\lambda)-(\Delta_{g_n}+\lambda_n)\right)v.
\end{equation}
We have that $\eta_n$ converges to $\eta$ and $(q_n,c_n,f_n)$ to $(q,c,f)$ in the smooth topology. Then, the right hand side of \eqref{one more} goes to $0$ in the smooth topology. This means that so does the left hand side. Observe now that, since $q_n$ converges to $q$, by \cref{schauder estimates} there exists a uniform positive constant $C$ such that
\begin{equation}\label{truly last}
    \|v_n-v\|_{H_{k+2}(M)}\le C\left(\|(\Delta_{g_n}+\lambda_n)(v_n-v)\|_{H_{k}(M)}+\|v_n-v\|_{L^2(M)}\right)
\end{equation}
for every sufficiently large $n$ and $k\geq\ell$. Since both terms on the right hand side of \eqref{truly last} go to $0$ for every $k$, $v_n$ converges to $v$ in the smooth topology. This proves continuity of the total map $D_q\tilde\Psi^{-1}$.
 
By construction, $D_q\tilde\Psi^{-1}$ is a linear inverse of $D_q\tilde\Psi$ for every $q$. We now want prove that $D\tilde\Psi^{-1}$ is smooth and tame in a neighborhood of $q_0$ in $\mathcal V$. By \cref{conti smooth}, it suffices to prove that it is tame. Tameness is a local property. Thus, it is enough to restrict to an open set in the codomain where $\|q\|_{0},|\langle u,v\rangle_{L^2_g}|$ and $\|(\Delta_g+\lambda)v+\eta u\|_{H_0}$ are bounded by a constant $C$. We first make estimates for $n=0$. By \cref{elliptic estimates}, we have
\begin{equation}\label{inductive estimate}
     \|v\|_{H_2}+|\eta|\leq C_0\left(\|(\Delta_g+\lambda)v+\eta u_0\|_{H_0}+\left|\langle u_0,v\rangle_{L^2_{g_0}}\right|\right).
\end{equation}
Fix a positive constant $\epsilon<\frac{1}{C_{0}}$. Then,
\begin{align}\label{uniform bound}
     \|(\Delta_g+\lambda)v+\eta u_0\|_{H_0}\leq& \|(\Delta_g+\lambda)v+\eta u\|_{H_0}+|\eta|\cdot\|u-u_0\|_{H_0}\nonumber
     \\ \leq& \|(\Delta_g+\lambda)v+\eta u\|_{H_0}+\epsilon|\eta| \nonumber
     \\ \leq& C+\epsilon|\eta|
\end{align}
for $(g,\lambda, u)$ in a small neighborhood of $q_0$. Thus, \eqref{inductive estimate} gives
\begin{equation}\label{zero estimate}
\|v\|_{H_2}+|\eta|\leq \tilde{C_0}\end{equation}
in the chosen open set. In particular, \eqref{uniform bound} gives a uniform bound on $\|(\Delta_g+\lambda)v+\dot\lambda u_0\|_{H_0}$. This gives estimates for $n=0$. Let now $n$ be strictly positive. By \cref{elliptic estimates}, we have
\begin{align}\label{main estimate}
\|v\|_{H_{n+2}}\leq C_n\left(1+\|(\Delta_g+\lambda)v+\dot\lambda u_0\|_{H_n}+\|(\Delta_g+\lambda)v+\eta u_0\|_{H_0}\|q\|_{n}+\|q\|_{n}\right),
\end{align}
where we have already applied $\left|\langle u_0,v\rangle_{L^2_{g_0}}\right|\leq C$. As in \eqref{uniform bound}, we obtain
\begin{align*}
    \|(\Delta_g+\lambda)v+\eta u_0\|_{H_n}\leq&\|(\Delta_g+\lambda)v+\eta u\|_{H_n}+|\eta|\cdot\|u-u_0\|_{H_n}
     \\ \leq& \|(\Delta_g+\lambda)v+\eta u\|_{H_n}+\tilde C\|q\|_{n}
\end{align*}
for a positive constant $\tilde C$. Thus, \eqref{main estimate} gives
\begin{equation*}
   \|v\|_{H_{n+2}}\leq \tilde C_n\left(1+\|(\Delta_g+\lambda)v+\eta u\|_{H_n}+\|q\|_{n}\right),
\end{equation*}
for a positive constant $\tilde C_n$ depending only on $n$. This proves tameness. By \cref{nash moser}, there is a smooth tame right inverse of $\Psi$ on $\mathcal V$. Therefore, $\mathcal V\cap \Psi^{-1}(U\times I\times\{0\})=\mathcal V \cap\mathcal Q'$ is the image of a smooth tame embedding, where $U$ is the projection of $\mathcal V$ on $\cG_+$ and $I\subseteq\R$ is a sufficiently small open interval containing $\frac 12 \|u_0\|^2_{L^2_{g_0}}$. Since this projection is an open map, $U$ is open. Therefore, $\mathcal Q'$ is a smooth tame Fre\'chet submanifold.
\end{proof}
From \cref{simple} it follows the smooth dependence of the eigenfunction on the metric. In particular, the following is an immediate consequence of the previous proof.
\begin{corollary}[Tame local dependence of simple eigenpairs]\label{tameness of eigen}
Let $q_0=(g_0,\lambda_0,u_0)\in Q'$, and assume that $\lambda_0$ is the
$k$--th eigenvalue of $g_0$. Then there exist an open neighborhood $U\subset \cG_+$ of $g_0$, an open interval $I\subseteq \R $, and an open neighborhood $\mathcal V\subset Q'$ of $q_0$ such that
\[
        \Psi|_\mathcal V : \mathcal V \longrightarrow U\times \mathbb{R}\times \{0\}
\]
is a smooth tame diffeomorphism onto its image. In particular, for every
$c\in I$, the map
\[
        \operatorname{eig}_{k,c}:U\longrightarrow Q',\qquad
        \operatorname{eig}_{k,c}(g):=(\Psi|_\mathcal V)^{-1}(g,c,0),
\]
is smooth and tame.
\end{corollary}
If $c=1/2$, then $\mathrm{eig}_{k,1/2}(g)$ gives a local smooth tame
choice of $L^2_g$-normalized eigenfunction.
We now describe the geometry of the connected components of $\mathcal Q'$. In particular, we show that the eigenpair space is divided in a family of subspaces corresponding to the index of the eigenvalue. Let $k$ be a positive integer. Denote by $\mathcal Q'_k$ the set of triples $(g,\lambda,u)$ in $\mathcal Q'$ where $\lambda$ is the $k$-th eigenvalue of $g$.
\begin{proposition}[Connected components of $\mathcal Q'$]\label{connected components'}
    For every positive integer $k$, the subspace $\mathcal Q_k'$ is union of connected components of $\mathcal Q'$. In particular, it is a tame Fr\'echet submanifold.
\end{proposition}
\begin{proof}
    By continuity of the $k$--th eigenvalue with respect to the metric, $\mathcal Q_k'$ is open in $\mathcal Q'$. This follows from \cite[Theorem 2.2]{Bando}. Since $\mathcal Q'$ is the disjoint union of $\{\mathcal Q'_k\}_{k\in\mathbb{N}}$, then $\mathcal Q'_k$ is also closed in $\mathcal Q'$. Therefore these sets are union of connected components of $\mathcal Q'$. In particular, they are tame Fr\'echet submanifolds of codimension $0$.
\end{proof}
The right inverse constructed in Section 2 is available on $\mathcal N$. We therefore restrict the eigenpair space to the following subset
\[\mathcal Q=\bigg\{(g,\lambda,u)\in\mathcal Q'\bigg|(g,u)\in \mathcal N\bigg \}.
\]
The following is an immediate consequence of \cref{Uhlenbeck}.
\begin{proposition}\label{tame}
    The set $\mathcal Q$ is open dense in $\mathcal Q'$. In particular, it is a smooth tame Fr\'echet submanifold of $\cG\times\R_{>0}\times C^\infty(M)$.
\end{proposition}
\begin{proof}
    If $\dim(M)> 2$, then this statement only claims that the generic eigenfunction of the Laplace--Beltrami operator is Morse, which follows from \cref{Uhlenbeck}. Assume now that $\dim(M)\geq 2$. Observe that if $(g,u)\in \mathcal N$ and $u$ is a nonconstant eigenfunction of $\Delta_g$, then the condition $\crit(u)\cap (\Delta_g)^{-1}(0)=\emptyset$ is equivalent to having $u$ transverse to $0$. Then the statement follows from \cref{Uhlenbeck} in the same way.
\end{proof}
We can now prove \cref{theorem A}. The rest of this work will be dedicated to discussing its topological consequences.
\begin{theorem}[Flexibility of eigenfunctions]\label{openness}
    The projection map $\pi\colon\mathcal Q\to C^{\infty}(M)$ is an open map.
\end{theorem}
\begin{proof}
    By \cref{simple}, $\mathcal Q$ is a tame Fr\'echet variety. Let $q_0$ in $\mathcal Q$ and $\mathcal{V}$ a neighborhood of $q_0$ as in \cref{laplacian use}. We want to show that $\pi$ is a smooth tame submersion with a smooth tame right inverse of $D\pi \colon T\mathcal V \to  C^\infty(M)$. Notice that a triple $(h,\eta,v)\in \cG\times \R\times C^\infty(M)$ belongs to $T_{q}\mathcal Q$, for $q=(g,\lambda, u)$, if and only if 
    \begin{equation} \label{expression for tangent}
        \dot \Delta_g[h](u)+\eta u+(\Delta_g+ \lambda)v=0.
    \end{equation}
   Let $R$ be a smooth tame right inverse of $\dot{\Delta}$ in $\mathcal V$ given by \cref{laplacian use}, i.e. $\dot \Delta_g[R(q,v)](u)=v$ for each $v\in C^\infty(M)$. Then setting $$\nu(q,v)=(R(q,-(\Delta_g+\lambda)v),0,v)$$ we have that $\nu(q,v)$ satisfies \eqref{expression for tangent}. Therefore
   \[\nu\colon \mathcal V \times C^\infty(M)\to T\mathcal V\]
   is a smooth tame  right inverse of $D\pi$ and it is linear on the second component.  We can conclude using \cref{nash moser}. The same proof applies to the projection on the eigenvalue and on both the eigenvalue and the eigenfunction.
\end{proof}
\begin{remark}[Banach setting]
    \cref{openness} might not hold if we set up the problem with $C^k$ metrics and $C^k$ eigenfunction. This is because the function $w_q(v)$ defined in the proof of \cref{sol in crit} might not be sufficiently regular. Indeed, the map $a_q$ in \eqref{simplified eq} would only be $C^{k-1}$, which means the solution constructed in \eqref{solution} might yield only a $C^{k-1}$ metric. The obstruction here is given by a loss of regularity.
\end{remark}
\begin{remark}[Image of the projection]
    As one would expect, the projection on $C^\infty(M)$ is not a surjective map. Its image is not even dense. Indeed, continuous maps with dense image send dense sets in dense sets. However, every eigenfunction of a Laplace--Beltrami operator has nonempty zero set. This is clearly not a dense property of smooth functions in general.
\end{remark}
\subsection{Generic properties of eigenfunctions}
For a smooth metric $g$ on $M$, denote by $[u_k(g)]$ its $k$-th eigenspace and by $\lambda_k(g)$ its $k$--th eigenvalue. We now aim to prove \cref{theorem b}. In other words, that every generic properties for smooth functions is also generic for eigenfunctions of the Laplace--Beltrami operator. To do so, it will be convenient to work with the projectivization of the eigenpair spaces.
\begin{definition}
    We define the projectivization of the eigenpair space of the Laplace--Beltrami operator as
    \[\tilde{\mathcal Q}=\{(g,\lambda,[u]):(g,\lambda,u)\in \mathcal Q\}\subseteq \cG_+\times\R_{>0}\times\mathbb PC^\infty.\]
    For $k\geq 1$, denote by $\tilde{\mathcal Q}_k$ the projectivization of $\mathcal Q_k$.
\end{definition}
\begin{remark}[Connected components of $\tilde{\mathcal Q}$]
    By \cref{connected components'}, $\tilde{\mathcal Q}_k$ is union of connected components of $\tilde{\mathcal Q}$.
\end{remark}
We are now ready to prove \cref{theorem b}.
\begin{theorem}[Generic properties of the $k$--th eigenfunction]\label{residual first}
    Let $M$ be a smooth compact connected manifold of dimension $\mathrm{dim}(M)\geq 2$ and $k$ be a positive integer. Let $D\subset \mathbb{P}C^{\infty}(M)$ and define
    \[\cG_+(k,D)=\{g\in \cG_+\colon \lambda_k(g)\mathrm{\ is\ simple}, [u_k(g)]\in D\}.\]
    The following hold
    \begin{enumerate}
        \item if $D$ is open, $\cG_+(k,D)$ is open.
        \item if $D$ is dense, $\cG_+(k,D)$ is dense.
        \item if $D$ is residual, $\cG_+(k,D)$ is residual.
    \end{enumerate}
\end{theorem}
\begin{proof}
    The first part of the statement is an immediate consequence of the continuity of the $k$--th eigenspace and $k$--th eigenvalue on the metric. This is given by \cref{tameness of eigen}. Let now $\bar{g}$ be a fixed smooth metric on $M$. Define $\mathcal Q_{\bar{g},k}=\{(g,\lambda,u)\in \mathcal Q_{k}:\|u\|_{L^{2}_{\bar g}}=1\}$. Consider the restriction and corestriction of the projection $\pi$
    \[\pi_{k}:\mathcal Q_{\bar{g},k}\to S_{L^{2}_{\bar g}}C^{\infty}(M).\]
    Since $\pi$ is an open map by \cref{openness}, $\pi_k$ is open as well. Indeed, if $f:X\to Y$ is an open continuous map and $Z\subseteq Y$, then $f|_{f^{-1}(Z)}:f^{-1}(Z)\to Z$ is also open continuous. Define the map
    \[p_k\colon \mathcal Q_{\bar g,k}\to \cG_+\times\R_{>0}\times\mathbb{P}C^{\infty}(M)\]
    as $p_k(g,\lambda,u)=(g,\lambda,p(u))$, where $p(u)$ is the vector span of $u$. Observe that $\tilde{\mathcal Q}_k=p_k(\mathcal Q_{\bar g,k})$. The map $p_{k}$ is a local homeomorphism with its image. Denote by $\tilde{\pi}_{k}:\tilde{\mathcal Q}_{k}\to\mathbb{P}C^{\infty}(M)$ the projection on the last component. We have the following commutative diagram.
\begin{center}\begin{tikzpicture}[>=stealth, node distance=2cm]
  \node (A) at (0,2) {$\mathcal Q_{\bar g,k}$};
  \node (B) at (0,0) {$\tilde {\mathcal Q}_k$};
  \node (C) at (4,2) {$S_{L^{2}_{\bar g}}C^{\infty}(M)$};
  \node (D) at (4,0) {$\mathbb{P}C^\infty(M)$};
  \draw[->] (A) --node[left]{$p_k$} (B);
  \draw[->] (A) --node[above]{$\pi_k$} (C) ;
  \draw[->] (B) --node[below]{$\tilde \pi_k$} (D);
  \draw[->] (C) --node[right]{$p$} (D) ;
\end{tikzpicture}\end{center}
    The maps $p_k$ and $p$ are local homeomorphism. Then $\tilde{\pi}_k$ is open if and only if $\pi_k$ is open. This is because being open is a local property. Thus, $\tilde \pi_k$ is open. Consider now the following diagram,
\begin{center}\begin{tikzpicture}[>=stealth, node distance=2cm]
  \node (A) at (2,2) {$\tilde {\mathcal Q}_k$};
  \node (B) at (0,0) {$\cG_+$};
  \node (C) at (4,0) {$\mathbb{P}C^\infty(M)$};
  \draw[->] (A) --node[left]{$P_k$} (B);
  \draw[->] (A) --node[above]{$\tilde \pi_k$} (C) ;
\end{tikzpicture}\end{center}
    where $P_k$ is the projection on the first component. The map $P_k$ is a homeomorphism with its image. Indeed, it is continuous and injective. To show that the inverse is continuous as well, we need to show that the $k$-th eigenvalue and eigenspace depend continuously on the metric. This follows from \cref{tameness of eigen}. Then
    \[\tilde\pi_k\circ P_{k}^{-1}:P_{k}(\tilde{\mathcal Q}_k)\subseteq\cG_+\to\mathbb{P}C^{\infty}(M)\]
    is an open continuous map. Observe that $P_{k}(\tilde{\pi}_{k}^{-1}(D))\subseteq \cG_+(k,D)$ by definition. Assume now that $D$ is dense. Standard topology gives that the preimage of a dense set under an open map is dense. Thus, $P_{k}(\tilde{\pi}_{k}^{-1}(D))$ is dense in $P_{k}(\tilde {\mathcal Q}_k)$. Since $P_{k}(\tilde {\mathcal Q}_k)$ is dense in $\cG_+$, the set $P_{k}(\tilde\pi_{k}^{-1}(D))$ is dense in $\cG_+$. This proves the second part of the statement. For the third part, choose open dense sets $O_j$ with $\bigcap_jO_j\subset D$. Each $\cG_+(k,O_j)$ is open dense by the first two parts, and their intersection is contained in $\cG_+(k,D)$. This proves residuality with the convention used here.
\end{proof}
We can use \cref{residual first} to prove an similar result for metrics of finite regularity. Recall that we denote by $\cG_+^\ell$ the set of $C^\ell$ Riemannian metrics on $M$. By \cite[Theorem 6.17]{Gilbarg}, eigenfunctions of Laplace--Beltrami operators associated to one such metric are in $C^{\ell}(M)$.
\begin{theorem}[Generic properties of the $k$--th eigenfunction for $C^\ell$-- metrics]\label{residual second}
    Let $M$ be a smooth compact connected manifold of dimension $\mathrm{dim}(M)=n\geq 2$ and $k$ be a positive integer and $\ell\geq 2 $ be a nonnegative integer. Let $D\subset \mathbb{P}C^{\ell}(M)$ and define
    \[\cG_+^\ell(k,D)=\{g\in \cG_+^\ell\colon \lambda_k(g)\mathrm{\ is\ simple}, [u_k(g)]\in D\}.\]
    The following hold
    \begin{enumerate}
        \item if $D$ is open, $\cG_+^\ell(k,D)$ is open.
        \item if $D$ is open and $i_\ell^{-1}(D)$ is dense in
       $\mathbb P C^\infty(M)$, then $\mathcal G_+^\ell(k,D)$ is open dense.
 \item Suppose $D$ contains $\bigcap_{j\geq1}D_j$, where every $D_j$
       is open in $\mathbb P C^\ell(M)$ and every
       $i_\ell^{-1}(D_j)$ is dense in $\mathbb P C^\infty(M)$.
       Then $\mathcal G_+^\ell(k,D)$ is residual.
    \end{enumerate}
\end{theorem}
\begin{proof}
    We first define the following subset of $\cG_+^k\times \R_{>0}\times \mathbb PC^{\ell}$
    \[\hat {\mathcal Q}_{k,C^\ell}=\{(g,\lambda,[u]):(\Delta_g+\lambda)u=0,\ m(g,\lambda)=1,(g,u)\text{ is nonresonant},\ \lambda\ \mathrm{is\ the\ }k-th\ \mathrm{eigenvalue}\}.\]
    Consider the following commutative diagram
\begin{center}\begin{tikzpicture}[>=stealth, node distance=2cm]
  \node (A) at (0,2) {$\mathcal {\tilde Q}_k$};
  \node (B) at (0,0) {$\hat {\mathcal Q}_{k,C^\ell}$};
  \node (C) at (4,2) {$\mathbb PC^\infty$};
  \node (D) at (4,0) {$\mathbb{P}C^\ell$};
  \draw[->] (A) --node[left]{$i_\ell$} (B);
  \draw[->] (A) --node[above]{$\tilde\pi$} (C) ;
  \draw[->] (B) --node[below]{$\tilde\pi_{C^\ell}$} (D);
  \draw[->] (C) --node[right]{$i_{\mathbb P_\ell}$} (D) ;
\end{tikzpicture}\end{center}
where the $\pi$ arrows are projections onto the third component and the $i$ arrows are inclusions. Observe that all these maps are continuous and the vertical ones are injections with dense image. Moreover, $\tilde{\pi}$ is open. This was shown in the proof of \cref{residual first}. Consider now a set $D\subseteq \mathbb{P}C^{\ell}$. By identifying $\mathcal {\tilde Q}_{k,C^\ell}$ the with the set of $C^\ell$ metrics, as we did in the proof of \cref{residual first}, we reduce to proving the statement for $\tilde\pi^{-1}_{C^\ell}(D)$. If $D$ is open, then $\tilde\pi_{C^\ell}^{-1}(D)$ is open by continuity of the horizontal map. This proves the first part of the statement. Assume now that $D$ is open dense. Since 
\begin{equation}\label{supset}
    \tilde\pi_{C^\ell}^{-1}(D)\supseteq i_{\ell}(i_{\ell}^{-1}(\tilde\pi_{C^{\ell}}^{-1}(D)))=i_\ell(\tilde\pi^{-1}(i_{\mathbb P_\ell}^{-1}(D))),
\end{equation}
where in the second equality we used the commutativity of the diagram. By hypothesis, the preimage $i_{\mathbb P_\ell}^{-1}(D)$ is open dense. Then, since $\tilde\pi$ is open, $\tilde\pi^{-1}(i_{\mathbb P_\ell}^{-1}(D))$ is also open dense. This is because the preimage of a dense set over an open map is still dense. Finally, since $i_{\ell}$ is a continuous map with dense image, $i_\ell(\tilde\pi^{-1}(i_{\mathbb P_\ell}^{-1}(D)))$ is dense in $\mathcal {\tilde Q}_{k,C^\ell}$. Then, by \eqref{supset}, $\tilde{\pi}_{C^\ell}^{-1}(D)$ is dense. This proves the second part of the statement. The third part of the statement immediately follow from the second after taking countable intersections.
\end{proof}
\begin{remark}
The second hypothesis in $(2)$ of \cref{residual second} cannot be dropped. Indeed, fix $0<\alpha<1$, $\ell\geq 2$ and a smooth reference metric $\bar g$. Define the following subsets of
$\mathbb P C^\ell(M)$
\[
 K_m=\left\{[u]:\frac{\|u\|_{C^{\ell,\alpha}}}{\|u\|_{L^2_{\bar g}}}
                    \leq m\right\},\qquad m\geq1.
\]
It is easy to see that $K_m$ is compact using Arzel\`a--Ascoli \cite[Theorem 3.25]{MR681859}. Moreover, an infinite-dimensional Banach manifold has no compact subset with nonempty interior. Thus
\[
 D=\mathbb P C^\ell(M)\setminus\bigcup_{m\geq1}K_m
\]
is residual. Nevertheless, elliptic regularity \cite[Theorem 6.17]{Gilbarg} implies that no
eigenfunction of any $C^\ell$ metric belongs to $D$. Thus $\mathcal G_+^\ell(k,D)=\emptyset$. This is why the compatibility assumption in \cref{residual second} is essential.
\end{remark}
Taking a countable intersection gives the following corollary.
\begin{theorem}[Generic properties of every eigenfunction]\label{residual}
    Let $M$ be a smooth compact connected manifold of dimension $\mathrm{dim}(M)\geq 2$. Let $D\subset \mathbb{P}C^{\infty}(M)$ be residual and define
    \[\cG_+(D)=\{g\in \cG_+\colon \lambda_k(g)\mathrm{\ is\ simple}\ \mathrm{and}\ [u_k(g)]\in D\mathrm{\ for\ every\ }k\geq 1\}.\]
    Then, $\cG_+(D)$ is residual.
\end{theorem}
The remaining extensions are stated in the smooth category. Their genericity assertions use the smooth openness theorem.
\subsection{Generic properties of metrics induced by an embedding}
A natural way of equipping $M$ with a Riemannian metric is to embed it into the Euclidean space and pull back the Euclidean metric. Let $\dim(M)=n$ and $N\geq 2n+1$ be a positive integer. Denote by $\mathrm{Emb}(M,\mathbb{R}^N)$ the space of smooth embeddings of $M$ into $(\mathbb R^N,g_N)$, where $g_N$ is the Euclidean metric on $\mathbb R^N$. This is an open dense subset of $C^{\infty}(M,\mathbb R^N)$ by the standard transversality theorem for embeddings. Given an embedding $f_0$ of $M$, we can consider the metric $f_0^\ast g_N$ on $M$ given by the pullback of the Euclidean metric. This induces a map between Fr\'echet manifolds
\begin{equation}\label{pullback}
    \mathcal P:\mathrm{Emb}(M,\mathbb{R}^N)\to \mathcal G^+
\end{equation}
defined as $\mathcal P(f)=f^*g_N$. This map was extensively studied by J. Nash in \cite{Nash} to prove his celebrated isometric embedding theorem.
\begin{theorem}[Nash's isometric embedding theorem {\cite[Theorem 3.1.1]{NashOverview1}}]\label{isom embedding}
    Let $3\leq k\leq \infty$, $n\geq 1$ and $N=\frac{n(3n+11)}{2}$. If $(M,g)$ is a compact $C^k$ Riemannian manifold of dimension $n$ with a $C^k$ metric, then there exists a $C^k$ isometric embedding $f:(M,g)\to\mathbb (R^N,g_N)$.
\end{theorem}
We refer to \cite{NashOverview2} and \cite{NashOverview1} for a survey of this result. A key ingredient in Nash's proof is the local openness of the metric-inducing map at free embeddings. We aim to use that approach to prove the following.
\begin{theorem}[Generic properties of induced metrics]\label{embedding residual}
    Let $M$ be a smooth compact connected manifold of dimension $\mathrm{dim}(M)=n\geq 2$. Let $N$ and $k$ be positive integers and assume $N\geq\frac{n(n+5)}{2}$. Let $D\subset \mathbb{P}C^{\infty}(M)$ and define
    \[E(k,D)=\{f\in \mathrm{Emb}(M,\mathbb R^N)\colon[u_k(f^\ast g_N)]\mathrm{\ is\ simple},[u_k(f^\ast g_N)]\in D\}.\]
    The following hold
    \begin{enumerate}
        \item if $D$ is open, $E(k,D)$ is open.
        \item if $D$ is dense, $E(k,D)$ is dense.
        \item if $D$ is residual, $E(k,D)$ is residual.
    \end{enumerate}
    Moreover, if $D$ is residual, then the set
    \[E(D)=\{f\in \mathrm{Emb}(M,\mathbb R^N)\colon \lambda_k(f^\ast g_N)\mathrm{\ is\ simple},[u_k(f^\ast g_N)]\in D\mathrm{\ for\ every\ }k\}\]
    is residual.
\end{theorem}
To prove \cref{embedding residual}, it is enough to show that the pullback map is open in an open dense subset of $\mathrm{Emb}(M,\mathbb R^N)$ under the hypothesis of the theorem. Since the preimage of an open dense set under an open continuous map is still open dense, the result then follows directly from \cref{residual first} and \cref{residual}. Following Nash's work, we define the following property.
\begin{definition}[Free embedding]
    Let $f_0$ be a smooth embedding of $M$ in $\mathbb R^N$. We say that $f_0$ is free if, for every $x\in M$ and every coordinate system $\{x_1,\ldots,x_n\}$ in a neighborhood of $x$, the following vectors are linearly independent
    \begin{equation}\label{free vectors}
        \frac{\partial f_0}{\partial x_i}(x), \frac{\partial^2 f_0}{\partial x_ix_j}(x)\quad 1\leq i\leq j\leq n.
    \end{equation}
\end{definition}
\begin{remark}
    Observe that freedom is an intrinsic condition on the $2$--jets of $f_0$. In particular, it does not depend on the choice of atlas on $M$. Since $M$ is compact, the set of free embeddings is open in the $C^2$ topology, and thus in the smooth one as well.
\end{remark}
No free embeddings can exist if $N<\frac{n(n+3)}{2}$, which is the number of vectors in \eqref{free vectors}. However, we need to ask for a stricter bound on the dimension of the target space. This stronger bound is needed because the proof of \cref{embedding residual} requires free embeddings to be dense. We refer to \cite{Ambra} for the following theorem although it should be attributed to J. Nash and M. Gromov.
\begin{proposition}[{\cite[Theorem N1, Theorem N2, Theorem IFT]{Ambra}}]\label{openness of embeddings}
    Let $M$ be a compact smooth manifold of dimension $n$ and let $N\geq\frac{n(n+5)}{2}$ be a positive integer. Then the set of free embeddings of $M$ in $\mathbb R^N$ is open dense in $\mathrm{Emb}(M,\mathbb R^N)$. Moreover, the pullback map $\mathcal P$ is open when restricted to the set of free embeddings.
\end{proposition}
\cref{openness of embeddings} only gives openness of the pullback over free embeddings. However, since this space is open dense, \cref{embedding residual} immediately follows.
\subsection{Schr\"odinger operators}\label{sect Schro}
So far we've restricted our study to the Laplace--Beltrami operators. In this section we consider a more general family of self--adjoint second order elliptic operators on $M$.
\begin{definition}[Schr\"odinger operator]
    Let $M$ be a smooth manifold. We say that a partial differential operator on $M$ is a Schr\"odinger operator if it is in the form
    \[\Delta_g+v,\]
    where $g$ is a smooth metric and $v$ is a smooth function on $M$
\end{definition}
Schr\"odinger operators are clearly selfadjoint and elliptic. Moreover, their spectrum is a nonincreasing sequence 
\[
-\lambda_1(g,v)\ge -\lambda_2(g,v)\ge \cdots \to -\infty.
\]
Thus, it makes sense to say which is the $k$--th eigenvalue of a Schr\"odinger operator. If $\lambda$ is a scalar, we will denote by $m(g,v,\lambda)=\dim(\mathrm{ker}(\Delta_g+v+\lambda))$. We will parametrize them naturally as\[(g,v)\to \Delta_g+v\] and study the generic properties of this family. In particular, we want to prove an analogue of \cref{openness} for Schr\"odinger operators. We first define the corresponding eigenpair space
\[
\mathcal S'=\{(g,v,\lambda,u):(\Delta_g+v+\lambda)u=0,m(g,v,\lambda)=1,u\neq 0\}\subseteq \cG_+\times C^\infty(M)\times\R\times C^\infty(M).
\]
The proof of \cref{simple} applies in the same way, giving that $S'$ is a smooth tame submanifold of the ambient space. Moreover, we define
\[\mathcal S=\{(g,v,\lambda,u)\in \mathcal S':(g,u)\in\mathcal N\}.\]
This is the subset of the Schr\"odinger eigenpair space where the differential of the Laplace--Beltrami operator is locally surjective with smooth tame right inverse. We want to show that $\mathcal S$ is open dense in $\mathcal S'$. In preparation for what follows, it will be convenient to define $S_{\ell}'$ to be the set of $(g,v,\lambda,u)$ in $\mathcal S'$ where $\lambda$ is the $\ell$--th eigenvalue of $\Delta_g+v$. Specifying the index of the eigenvalue this way allows us to identify the generic operator with its eigenspace. As we did in \cref{connected components'}, one can show that $\mathcal S_\ell'$ is union of connected components of $\mathcal S'$. In particular, it is open. We will prove that $\mathcal S_\ell$ is open dense in $\mathcal S'_\ell$ for every $\ell$. The following result will be useful in this regard. We only state a simplified version, which will be enough for our purposes. Recall that a Baire space is a topological space where every countable intersection of open dense sets is dense. Since every space we work with is Baire, we will not discuss this property further.
\begin{theorem}[Kuratowski-Ulam, {\cite[Theorem 3.5.16]{Srivastava}}]\label{Kurat}
    Let $X$ and $Y$ be second countable Baire spaces and $A\subseteq X\times Y$ be open. Then the following are equivalent:
    \begin{enumerate}
        \item $A$ is dense.
        \item The set $\{x\in X:A_x\ is\ residual\}$ is residual, where
        \[A_x=\{y\in Y:(x,y)\in A\}.\]
    \end{enumerate}
\end{theorem}
\cref{Kurat} allows to reduce global density to density in the fibers. To prove the following result, we define the finite-regularity eigenpair space of the Schr\"odinger operators as follows.
\begin{definition}
    Let $k\geq 2$ be a positive integer. We define the $C^k$ eigenpair space of Schr\"odinger operators as
    \[\mathcal S'_{C^k}=\{(g,v,\lambda,u):(\Delta_g+v+\lambda)u=0,m(g,v,\lambda)=1,u\neq 0\}\subseteq \cG^k_+\times C^{k-1}(M)\times\R\times H_{k}(M).\]
\end{definition}
\begin{remark}
    Recall that the eigenfunctions of an elliptic operator with $C^{k-1}$ coefficients are in $H_{k}$ by \cref{schauder estimates}. This justifies the spaces chosen in the definition above. Also observe that $\mathcal S'_{C^k}$ is a Banach manifold. This can be proven as in \cref{simple}.
\end{remark}
\begin{proposition}\label{morse 2}
    Let $M$ be a smooth compact connected manifold with $\dim(M)\geq 2$. For the generic $q\in \mathcal S'$, $u$ has nonsingular zero set, is Morse, and $(\Delta_gu)^{-1}(0)\cap \crit(u)=\emptyset$. Thus, the set $\mathcal S$ is open dense in $\mathcal S'$. In particular, it is a tame Fr\'echet submanifold.
\end{proposition}
\begin{proof}
    Let $\ell$ be a positive integer and $\dim(M)=n$. We will prove that the set of $q$ as in the statement is open dense in $\mathcal S'_\ell$. Let $k\geq n+3$ be a positive integer. Consider the space $\cG_+\times C^{k-1}(M)$ of Schr\"odinger operators given by a smooth metric and a $C^k(M)$ function. Denote by $A_{\ell,C^k}\subseteq \cG_+\times C^k(M)$ the set of Schr\"odinger operators with simple $\ell$--th eigenvalue and the corresponding eigenfunction Morse and transverse to $0$. This set is open. The fact that having simple $\ell$-th eigenvalue is an open condition is proven exactly like it is for the Laplace-Beltrami operator. By \cref{fiber uhlenbeck}, every fiber $A_{\ell,C^k,g}\subseteq C^k(M)$ is residual. Then, by \cref{Kurat}, the set $A_{\ell,C^k}\subseteq \cG_+\times C^k(M)$ is also residual. It is also open by compactness of $M$. Since this holds for every $k$ sufficiently large, using \cref{lemma Q_k residual} we obtain that
    \[A_\ell=A_{\ell,C^k}\cap (\cG_+\times C^{\infty}(M))\subseteq \cG_+\times C^{\infty}(M)\]
    is also open dense. Consider now the continuous map
    \[
    P_\ell:\mathcal S'_\ell\to \cG_+\times C^{\infty}(M)
    \]
    defined as $P_\ell(g,v,\lambda,u)=(g,v)$. Since the space of Schr\"odinger operators with simple $\ell$--th eigenvalue is open and $P_\ell$ is a projection, this map is open. In particular, the preimage of an open dense set will still be open dense. Thus, $P_\ell^{-1}(A_\ell)$ is open dense.
    
    It remains only to prove that the generic element in $S'$ avoids the third condition of the statement. Consider the map
    \[F:\mathcal S'_{C^k}\times M\to T^*M\oplus \R\]
    defined as $F(q,x)=(x,d_xu,\Delta_gu(x))$. We want to show that $F$ is transverse to the set
    \[Z=\{(x,0,0):x\in M\}\]
    and then apply \cref{transversality theorem}. Assume then that $(g,v,\lambda,u)\in S'_{C^k}$ and $x\in M$ are such that $F((v,\lambda,u),x)\in Z$. By arguing as in the first part of the proof, we can show that for every $(\tilde g,\tilde v,\tilde \lambda,\tilde u)$ in an open dense subset of $\mathcal S'_{C^k}$, $\tilde u$ is transverse to $0$. To complete the proof, it is enough to prove that the restriction of $F$ to this set is transverse to $Z$. Thus, we can assume that there exists an open neighborhood $V$ of $x$ where $u$ is never vanishing. Observe that the tangent to $\mathcal S'_{C^k}$ in $(g,v,\lambda,u)$ is defined by the equation
    \[\dot\Delta_{g}[\dot g](u)+(\dot v+\dot\lambda )u+(\Delta_g+v+\lambda)\dot u=0.\]
    To prove submersivity of $F$, it is enough to consider variations with $\dot g|_V=0$ and $\dot\lambda=0$. Indeed, let $w$ be any smooth function compactly supported in $V$. Consider the function
    \[\dot v(w)=-\frac{(\Delta_g+v+\lambda)w}{u}.\]
    By construction, the element $(0,\dot v(w),0,w)$ verifies the equation characterizing the tangent of $\mathcal S'_{C^k}$. To prove that it belongs to the tangent, we need only to prove that $\dot v(w)\in C^{k-1}(M)$. Recall that throughout this work we always assume $k$ to be sufficiently large such that $u\in C^2(M)$. Then, since $u$ is a $C^2$ eigenfunction of an elliptic operator with $C^{k-1}$ coefficients, we have $u\in C^k$. This is a consequence of \cite[Theorem 6.17]{Gilbarg}. Moreover, $(\Delta_g+v+\lambda)w$ is also a $C^{k-1}$ function. This is because $w$ is smooth and $\Delta_g+v+\lambda$ has $C^{k-1}$ coefficients. This proves that $(0,\dot v(w),0,w)\in T_q \mathcal S_{C^k}'$. Then $F$ is a submersion in $(q,x)$, since we can always consider a smooth function $w$ compactly supported in $V$ with prescribed differential on $\Delta_gw$ in $x$. Since $\codim (Z)=n+1>\dim (M)$, by \cref{transversality theorem} we have that $\mathcal S'_{C^k}\cap \{q:F(q,M)\cap Z=\emptyset\}$ residual in $\mathcal S'_{C^k}$. It is also open by compactness of $M$. Since this holds for every $k$ sufficiently large, by \cref{lemma Q_k residual} we obtain that $\left(\mathcal S'_{C^k}\cap \{q:F(q,M)\cap Z=\emptyset\}\right)\cap \mathcal S'$ is open dense in $S'$. It is exactly the set of elements in $\mathcal S'$ avoiding the third condition in the statement. This concludes the proof.
\end{proof}
\begin{remark}
    Observe than the generic Schr\"odinger operator has simple eigenvalue by \cref{fiber uhlenbeck} even in dimension $1$. This is very different from the behaviour of the Laplace--Beltrami operators, whose nonzero eigenvalues always have multiplicity $2$. 
\end{remark}
We can finally state the main result of this section. The proof is an immediate consequence of \cref{laplacian use} and \cref{morse 2}. It essentially coincides with the one of \cref{openness}.
\begin{theorem}[Flexibility of eigenfunctions for Schr\"odinger operators]\label{schrodinger openness}
    Let $M$ be a smooth compact connected manifold $\dim(M)\geq 2$. Then, the projection on the eigenfunction
    \[\pi_\mathcal S:\mathcal S\to C^\infty(M)\]
    is open. In particular, it is open in an open dense subset of $\mathcal S'$.
\end{theorem}
Analogues of \cref{residual first} and \cref{residual} follow immediately with the same proofs. In the same way, one can also prove generic properties of Schrodinger operators induced by embeddings. We denote by $\lambda_k(f^\ast g,v)$ the $k$--th eigenvalue of the operator $\Delta_g+v$.
\begin{theorem}[Generic properties of induced Schr\"odinger operators]
    Let $M$ be a smooth compact connected manifold with $\dim(M)=n$ and assume $n\geq 2$. Let $N$ and $k$ be positive integers and assume $N\geq\frac{n(n+5)}{2}$. Let $D\subset \mathbb{P}C^{\infty}(M)$ and define
    \[\mathcal S(k,D)=\{(f,v)\in \mathrm{Emb}(M,\mathbb R^N)\times C^\infty(M)\colon\lambda_k(f^\ast g_N,v)\mathrm{\ is\ simple},[u_k(f^\ast g_N,v)]\in D\}.\]
    The following hold
    \begin{enumerate}
        \item if $D$ is open, $\mathcal S(k,D)$ is open.
        \item if $D$ is dense, $\mathcal S(k,D)$ is dense.
        \item if $D$ is residual, $\mathcal S(k,D)$ is residual.
    \end{enumerate}
    Moreover, if $D$ is residual, then the set
    \[\mathcal S(D)=\{(f,v)\in \mathrm{Emb}(M,\mathbb R^N)\times C^\infty(M)\colon\lambda_k(f^\ast g_N,v)\mathrm{\ is\ simple},[u_k(f^\ast g_N,v)]\in D\mathrm{\ for\ every\ }k\}\]
    is residual.
\end{theorem}
\begin{remark}
    One could also consider different families of second order self--adjoint elliptic operators and apply the methods used in this paper. An interesting one is given by
    \[\{-4\Delta_g+\mathrm{Scal}_g:g\in\cG_+\},\]
    where $\mathrm{Scal}_g$ is the scalar curvature associated to the metric $g$. This gives an important class of Schr\"odinger operators, which are studied in \cite{perelman2002entropy}. However, since the function involved in the operator depends on the metric, one cannot expect to prove openness of the projection by simply applying \cref{openness}.
\end{remark}
\section{Appendix}
\subsection{Proof of \cref{lem:metric variation of laplacian}}
Consider a curve of metrics $g_t=g+th$. We make computations in coordinates. Thus, we can write
    \begin{equation}\label{coordinate laplacian}
        \Delta_{g_t}u=g_t^{ij}(\partial_i\partial_ju-\Gamma^{k}_{ij}(t)\partial_ku),
    \end{equation}
    where $\Gamma^{k}_{ij}(t)$ are the Christoffel symbols w.r.t. the metric $g_t$. This expression is obtained directly from the definition of $\Delta$. We have
    \[\frac{d}{dt}g_t^{ij}\Big|_{t=0}=-g^{ik}g^{j\ell}h_{k\ell}=-h^{ij}.\]
    This gives the first summand of  \eqref{differential of Lapl}. We now differentiate the Christoffel symbols. Recall
    \[\Gamma_{ij}^k(t)=\frac 12 g^{k\ell}_t(\partial_ig_{t,j\ell}+\partial_jg_{t,i\ell}-\partial_\ell g_{t,ij}).\]
    We will omit the time when $t=0$ to simplify notation. Differentiating gives
    \begin{align*}\dot\Gamma_{ij}^k&=-\frac 12 h^{k\ell}(\partial_ig_{j\ell}+\partial_jg_{i\ell}-\partial_\ell g_{ij})+\frac 12 g^{k\ell}(\partial_ih_{j\ell}+\partial_jh_{i\ell}-\partial_\ell h_{ij})\nonumber\\
    &=-h^{k\ell}g_{a\ell}\Gamma^a_{ij}+\frac 12 g^{k\ell}(\partial_ih_{j\ell}+\partial_jh_{i\ell}-\partial_\ell h_{ij}).\end{align*}
    We want to simplify the second term of the summand. We have 
    \[\partial_ih_{j\ell}+\partial_jh_{i\ell}-\partial_\ell h_{ij}=\nabla_ih_{j\ell}+\nabla_jh_{i\ell}-\nabla_\ell h_{ij}+C_{i,j,\ell},\]
    where $C_{i,j,\ell}$ is the correction given by the Christoffel symbols. In particular,
    \[C_{i,j,\ell}=\Gamma_{ij}^ah_{a\ell}+\Gamma_{i\ell}^ah_{aj}+\Gamma_{j\ell}^ah_{ai}+\Gamma_{ij}^ah_{a\ell}-\Gamma_{i\ell}^ah_{aj}-\Gamma_{j\ell}^ah_{ai}\nonumber=2\Gamma_{ij}^ah_{a\ell}.\]
    Moreover, observe that
    \[h^{k\ell}g_{a\ell}=g^{ik}g^{j\ell}h_{ij}g_{a\ell}=g^{ik}h_{ia}=g^{k\ell}h_{a\ell}.\]
    Thus, we have obtained
    \begin{align*}
    \dot\Gamma_{ij}^k&=\Gamma^a_{ij}(g^{k\ell}h_{a\ell}-h^{k\ell}g_{a\ell})+\frac 12 g^{k\ell}(\nabla_ih_{j\ell}+\nabla_jh_{i\ell}-\nabla_\ell h_{ij})\nonumber\\&=\frac 12 g^{k\ell}(\nabla_ih_{j\ell}+\nabla_jh_{i\ell}-\nabla_\ell h_{ij}).\end{align*}
    Now we want to compute $-g^{ij}\dot\Gamma^k_{ij}\partial_ku$. We have
    \begin{align*}
        g^{ij}\dot\Gamma^k_{ij}\partial_ku&=\frac 12g^{ij} g^{k\ell}(\nabla_ih_{j\ell}+\nabla_jh_{i\ell}-\nabla_\ell h_{ij})\partial_ku\\
        &=\frac 12g^{k\ell}(g^{ij}\nabla_ih_{j\ell}+g^{ji}\nabla_jh_{i\ell}-\nabla_\ell \tr_gh)\partial_ku\\
        &=g^{k\ell}(g^{ji}\nabla_jh_{i\ell}-\frac 12\nabla_\ell \tr_gh)\partial_ku\\
        &=g^{k\ell}(\mathrm{div}_g(h)_\ell-\frac 12\nabla_\ell \tr_gh)\partial_ku\\
        &=\langle\mathrm{div}_g(h)^\sharp-\frac 12\nabla \tr_gh,\nabla u\rangle_g.
    \end{align*}
    This gives \eqref{differential of Lapl}.
    
    Now we compute the formula with respect to a conformal variation. Let $\psi\in C^\infty(M)$. Consider the conformal variation $h=\psi g$. Then we have $-\langle \psi g,\nabla^2 u\rangle_g=-\psi \Delta_g u$. Moreover, 
    \[
    \frac12\nabla\tr_g \psi g=\frac n2 \nabla \psi\quad \text{ and}\quad -\dive_g (\psi g)^\sharp=-\nabla \psi,
    \] which implies equation \cref{conformal equation}. For a traceless variation $h$, we just observe the Leibniz rule
    $$\dive_g(h(\nabla u,\cdot))=\langle \dive_gh^\sharp,\nabla u\rangle_g+\langle h, \nabla^2 u\rangle_g.$$
    Substituting this into equation \eqref{differential of Lapl} gives \eqref{traceless equation}.

    We now prove \eqref{final 3.8}. By substituting $h=wF$ into \eqref{differential of Lapl} we get
\begin{align}\label{f laplacian}
    \dot \Delta_g[wF](u)=&-\langle h,\nabla^2 u\rangle_g+\langle \frac 1 2\nabla(\tr_g h)-\dive_gh^\sharp, \nabla u \rangle_g\nonumber\\
    =&-\langle F,\nabla^2 u \rangle_g w+\langle \frac 1 2\nabla(\tr_g Fw)-\dive_g (Fw)^\sharp, \nabla u\rangle_g.
\end{align}
We have
\begin{align}\label{riemann 1}
    \langle\nabla\tr_g Fw, \nabla u\rangle_g&=\langle\nabla(\tr_g F)w+\tr_g F\cdot\nabla w, \nabla u\rangle_g\nonumber\\
    &=\langle\nabla\tr_g F,\nabla u\rangle_gw+\langle\tr_g F\nabla u, \nabla w\rangle_g
\end{align}
and
\begin{align}\label{riemann 2}
    \langle\dive_g (wF)^\sharp, \nabla u\rangle_g&=\langle F(\nabla w,\cdot)^\sharp+(\dive_g F)w^\sharp, \nabla u\rangle_g\nonumber\\
    &=\langle F(\nabla w,\cdot)^\sharp,\nabla u\rangle_g+\langle\dive_g F^\sharp, \nabla u\rangle_gw.
\end{align}
Moreover, we have
\begin{align}\label{riemann 3}
    \langle F(\nabla w,\cdot)^\sharp,\nabla u\rangle_g=F(\nabla w,\nabla u)=F(\nabla u,\nabla w)=\langle F(\nabla u,\cdot)^\sharp,\nabla w\rangle_g.
\end{align}
Substituting \eqref{riemann 1}, \eqref{riemann 2} and \eqref{riemann 3} in \eqref{f laplacian} concludes the proof.\qed
\subsection{Proof of \cref{lyapunov estimates}}
Let $A\coloneqq JX_0(x_0)$. Then every eigenvalue of $A+\mu$ has negative real part. Thus, the matrix
    \[P=\int_0^{+\infty}e^{(A^T+\mu\mathrm{Id})t}e^{(A+\mu\mathrm{Id})t}dt\]
    is well defined, symmetric and positive definite. In particular, it satisfies the Lyapunov equation
    \[(A^T+\mu\mathrm{Id})P+P(A+\mu\mathrm{Id})=-\mathrm{Id},\]
    which can be written as
    \begin{equation}
        A^TP+PA=-\mathrm{Id}-2\mu P
    \end{equation}
    Refer to the proof of \cite[Theorem 4.6]{Khalil} for the details. Let $\epsilon>0$ be small enough such that
    \begin{equation}\label{epsilon}2\|P\|\epsilon<\frac12,\end{equation}
    where $\|P\|$ is the operator norm of the matrix $P$ associated to the standard Euclidean norm. Let $x\in\mathbb R^n$ and $B$ be a matrix such that $\|B-A\|<\epsilon$. Then,
    \begin{align}\label{small matrix diff}
        x^T(B^TP+PB)x&\leq x^T((B-A)^TP+P(B-A))x+x^T(A^TP+PA)x\nonumber\\&\leq 2\epsilon\|P\|\|x\|^2+ x^T(-\mathrm{Id}-2\mu P)x\nonumber\\&\leq-\frac{\|x\|}{2}^2-2\mu x^TPx\nonumber\\&\leq -2\mu x^TPx,
    \end{align}
    where we have used \eqref{epsilon} in the third inequality.

    We now construct a positively invariant neighborhood for $X_0$. Consider the function
    \[V_{0}(x)=(x-x_0)^TP(x-x_0).\]
    Let $r>0$ be a scalar. Since $P$ is positive definite, the ellipsoid
    $U=V_{0}^{-1}([0,r))$ centered in $x_0$ has compact closure and is convex. Moreover, we can pick $r$ small enough such that
    \[\|JX_0(x)-A\|\leq \frac{\epsilon}4\]
    for every $x\in \overline U$ and $X_0(x)\neq 0$ for every $x\in \overline U\setminus \{x_0\}$. Now we want to prove
    \begin{equation}\label{pointing inward}
        DV_{0}(x)[X_0(x)]\leq -2\mu V_{0}(x)<0
    \end{equation}
    for every $x\in\overline U\setminus \{x_0\}$. This will immediately imply that $U$ is a $\Phi_t^0$--positively invariant neighborhood. Proving the inequality in the middle is the only nontrivial step. We have
    \begin{align}
    X_0(x)&=X_0(x)-X_0(x_0)\nonumber\\&
    =\int_0^1\frac{d}{ds}X_0(x_0+s(x-x_0))ds\nonumber
    \\&=\int_0^1JX_0(x_0+s(x-x_0))ds(x-x_0)=B_x(x-x_0),
    \end{align}
    where $B_x=\int_0^1DX_0(x_0+s(x-x_0))ds$ is considered as an endomorphism of $\R^n$. This gives    
\begin{align}\label{gradient estimate}
    DV_{0}(x)[X_0(x)]&=(B_x(x-x_0))^TP(x-x_0)+(x-x_0)^TPB_x(x-x_0)\nonumber\\
    &=(x-x_0)^T(B_x^TP+PB_x)(x-x_0)
\end{align}
    By convexity of $U$, we have
    \begin{equation}\label{special estimate}
        \|B_x-A\|\leq\int_0^1\|JX_0(x_0+s(x-x_0))-A\|ds\leq \frac\epsilon4<\epsilon
    \end{equation}
    Therefore, we can substitute \eqref{small matrix diff} into \eqref{gradient estimate} and obtain that for every $x\in \overline U\setminus \{x_0\}$
    \begin{equation}\label{invariant neigh}
        DV_{0}(x)[X_0(x)]\leq -2\mu (x-x_0)^TP(x-x_0)=-2\mu V_{0}(x)<0.
    \end{equation}
    This means exactly that $U$ is $\Phi_t^0$--positively invariant. Indeed, assume by contradiction that there exists $x\in U$ and $t\geq 0$ such that $V_0(\Phi_0^t(x))>r$. Then there exists $t_0\geq 0$ s.t. $V_0(\Phi_0^{t_0}(x))$ is increasing in $t$, which contradicts \eqref{invariant neigh}. Finally, since $\overline U$ is compact and positively invariant, \cite[Theorem 3.3]{Khalil} implies that $\Phi_t^0(x)$ is defined for every $t\geq 0$ and $x\in U$.

    Now we show that this holds also for every vector field sufficiently $C^1$--close to $X_0$. Let $\delta$ be as follows. For every vector field $X_1$ such that $\|X_0-X_1\|_{C^1}<\delta$,
    \begin{equation}\label{last one}
        \|JX_1(x)-A\|\leq \frac{\epsilon}2
    \end{equation}
    for every $x\in \overline U$. Moreover, by \eqref{invariant neigh}, we can assume that
    \[DV_{0}(x)[X_1(x)]<0\]
    for every $x\in\partial U$ and $X_1$ as above. Such a $\delta$ exists by compactness of $\overline U$. Finally, we pick $\delta$ small enough such that every such $X_1$ has a critical point $x_1$ in $\overline U$. Since $X_0$ is nondegenerate in $0$ and $X_1$ is $C_1$--close to $X_0$, we can assume that $x_1$ is the only critical point of $X_1$ in $\overline U$. As before, this implies that $U$ is a positively invariant neighborhood for every such $X_1$ and that the flow of such vector fields starting at a point in $U$ is defined for every positive $t$. Fix one such $X_1$ and let $x_1$ be its equilibrium in $U$. We now prove the first part of \eqref{liapunov partial derivative}. Consider the quadratic form
    \[V_1(x)=(x-x_1)^TP(x-x_1).\]
    We want to prove that
    \begin{equation}\label{almost Gronwall}
        \frac{d}{dt}V_1(\Phi_t^1(x))=DV_1(\Phi_t^1(x))[X_1(\Phi^1_t(x))]\leq -2\mu V_1(\Phi_t^1(x))<0
    \end{equation}
    for every $x\in U\setminus \{x_1\}$. This is obtained exactly as we did to prove \eqref{invariant neigh}, but this time we need to use that $\Phi_t^1(x)$ belongs to $U$ to have the estimate \eqref{special estimate}. Applying Gr\"onwall lemma to \eqref{almost Gronwall} gives
    \[V_1(\Phi_t^1(x))\leq V_1(x)e^{-2\mu t}.\]
By definition of $V_1$, we have
\[    \lambda_\mathrm{min}(P)\|\Phi_t^1(x)-x_1\|^2\leq V_1(\Phi_t^1(x))\leq V_1(x)e^{-2\mu t}\leq\lambda_\mathrm{max}(P)\|x_1-x\|^2e^{-2\mu t}.\]
Thus, we have obtained
\begin{equation}\label{unique}
\|\Phi_t^1(x)-x_1\|\leq \sqrt{\frac{\lambda_\mathrm{max}(P)}{\lambda_\mathrm{min}(P)}}e^{-\mu t}\|x-x_1\|\leq c_0e^{-\mu t},\end{equation}
where $c_0=\sqrt{\frac{\lambda_\mathrm{max}(P)}{\lambda_\mathrm{min}(P)}}\mathrm{diam}(U)$ depends only on $X_0$ through $P$ and $U$.

It remains only to prove \eqref{liapunov partial derivative} for the flow of $X_1$. We prove it by induction. To simplify notation, we will write $X=X_1$ and denote its flow by simply $\Phi_t$. Differentiating the equation characterizing the flow of $X$, we obtain
\begin{equation}\label{m=1}
    D\dot\Phi_t(x)=D(X(\Phi_t(x)))=JX(\Phi_t(x))D\Phi_t(x).
\end{equation}
Observe that \eqref{m=1} is a linear system with associated matrix $JX(\Phi_t(x))$ and unknown $D\Phi_t(x)$. Since \eqref{last one} holds by construction for \eqref{m=1}, the Lyapunov estimates we gave for $m=0$ still hold. Explicitly, let $v\in \R^n$ be a vector of norm $1$ and define $Y_{v,t}(x)=D\Phi_t(x)v$. Then we have
\begin{equation}\label{m=1 Y}
    \dot Y_{v,t}(x)=JX(\Phi_t(x))Y_{v,t}(x).
\end{equation}
Proceeding as in the first part of the proof with the Lyapunov function \[V_{v,t}(x)=Y_{v,t}(x)^TPY_{v,t}(x)\]
and taking the supremum over $v$ gives
\begin{equation}\label{to copy}
    \|D\Phi_t(x)\|\leq \|D\Phi_0(x)\|c_1e^{-\mu t}=c_1e^{-\mu t},
\end{equation}
where the second inequality is given by $\Phi_0(x)=x$. This implies the estimate in \eqref{liapunov partial derivative}.

We now consider higher order derivatives. Let $m\geq 2$ and assume that the last part of the statement has been proven for every smaller integer. For simplicity, we will denote by $c_m$ every constant which can be made dependent only on $m$, $\delta$ and $X_0$. We have
\begin{equation}\label{main lyapunov estimate}D^m\dot\Phi_t(x)=D^m(X(\Phi_t(x)))=JX(\Phi_t(x))D^m\Phi_t(x)+F_{m,t}(x),
\end{equation}
where $F_{m,t}(x)$ is a finite sum of terms in the form
\[D^kX(\Phi_t(x))[D^{j_1}\Phi_t(x),\ldots,D^{j_k}\Phi_t(x)]\]
with $2\leq k\leq m$, $j_i\geq 1$ and $j_1+\ldots+ j_k=m$. The coefficients and number of these terms depend only on $m$. This means that
\begin{equation}\label{inequality Fm}
    \|F_{m,t}\|\leq c_m\sum_{k=2}^m\|D^kX\|\sum_{j_1+\ldots+j_k=m}\prod_{i=1}^k\|D^{j_i}\Phi_t\|.
\end{equation}
By induction, we have
\begin{equation}\label{inequality 1}
\|D^{j_i}\Phi_t\|\leq c_{j_i}e^{-\mu t}(1+\|X\|_{j_i})\leq c_me^{-\mu t}(1+\|X\|_{j_i}).
\end{equation}
Applying \eqref{inequality 1} to \eqref{inequality Fm} gives
\begin{align}\label{another one}
\|F_{m,t}\|&\leq c_m\sum_{k=2}^m\|X\|_k\sum_{j_1+\ldots+j_k=m}\prod_{i=1}^ke^{-\mu t}(1+\|X\|_{j_i})\nonumber\\&
\leq c_me^{-2\mu t}\sum_{k=2}^m\|X\|_k\sum_{j_1+\ldots+j_k=m}\prod_{i=1}^k(1+\|X\|_{j_i})\nonumber\\&
\leq c_me^{-2\mu t}\sum_{k=2}^m\sum_{j_1+\ldots+j_k=m}\left(\|X\|_k+\sum_{\ell=1}^k\sum_{f_1\neq\ldots\neq f_\ell}\|X\|_k\|X\|_{j_{f_1}}\cdots\|X\|_{j_{f_\ell}}\right),\end{align}
where all the combinatorial coefficients have been incorporated in the constant $c_m$. For $1\leq j\leq m$ and $m\geq 2$, \cref{interpolation} gives
\begin{equation*}
    \|X\|_{j}\leq c_m\|X\|_1^{\frac{m-j}{m-1}}\cdot\|X\|_m^{\frac{j-1}{m-1}}.
\end{equation*}
This gives
\begin{align}\label{lyapunov interpolation}
    \|X\|_k\|X\|_{j_{f_1}}\cdots\|X\|_{j_{f_\ell}}&\leq c_m \|X\|_1^{\frac{(\ell+1) m-(k+\sum j_f)}{m-1}}\cdot\|X\|_m^{\frac{k-1-\ell+\sum j_f}{m-1}}\nonumber\\&\leq c_m \|X\|_m^{\frac{k-1-\ell+\sum j_f}{m-1}}
\end{align}
Since all the $j$ indices appearing in \eqref{another one} are positive, we have
\[\sum_{i=1}^\ell j_{f_i}\leq \sum_{i=1}^k j_i-(k-\ell)=m-(k-\ell).\]
In particular, $0\leq \frac{k-1-\ell+\sum j_f}{m-1}\leq 1$. Then,
\begin{equation}\label{lyapunov interpolation 2}
    \|X\|_m^{\frac{k-1-\ell+\sum j_f}{m-1}}\leq 1+\|X\|_m.
\end{equation}
Substituting \eqref{lyapunov interpolation} and \eqref{lyapunov interpolation 2} in \eqref{another one} gives
\begin{align}\label{Fm}
    \|F_{m,t}\|&\leq c_me^{-2\mu t}\sum_{k=2}^m\sum_{j_1+\ldots+j_k=m}\left(\|X\|_k+\sum_{\ell=1}^k\sum_{f_1\neq\ldots\neq f_\ell}(1+\|X\|_m)\right)\nonumber\\
    &\leq c_me^{-2\mu t}(1+\|X\|_m)
\end{align}
where we estimated $\|X\|_k\leq \|X\|_m$ and incorporated all the combinatorial constants in $c_m$.

It remains only to estimate the first summand in the right hand side of \eqref{main lyapunov estimate}. We have
\[\|JX(\Phi_t)D^m\Phi_t\|\leq\|JX(\Phi_t)\|\cdot\|D^m\Phi_t\|\leq c_m\|D^m\Phi_t\|.\]
Therefore, it is enough to make estimates on $D^m\Phi_t$. To do so, we'll write it explicitly as a solution of \eqref{main lyapunov estimate}. Define the evolution operator $S(t,s,x)$ as the solution of 
\begin{equation}\label{evolution equation}
\dot S(t,s,x)=JX(\Phi_t(x))S(t,s,x),\quad S(s,s,x)=\mathrm{Id}_n\forall s\in\R_{\geq 0}, x\in \R^n.
\end{equation}
By the standard theory of ODEs, the solution of \eqref{evolution equation} exists and is unique given initial conditions. Moreover, \eqref{last one} holds for \eqref{evolution equation}. Thus, we can obtain the following Lyapunov estimates for $S(t,s,x)$ in the same way as we did for \eqref{to copy}
\begin{equation}\label{s estimate}
    \|S(t,s,x)\|\leq c_m e^{-\mu (t-s)}\|S(s,s,x)\|=c_m e^{-\mu (t-s)},
\end{equation}
where the first equality is an immediate consequence of \eqref{evolution equation}. An easy computation shows that the function
\begin{equation}\label{expression for derivative of flow}
    f_t(x)=\int_0^tS(t,s,x)F_{m,s}(x)ds
\end{equation}
is a solution to \eqref{main lyapunov estimate}. In particular, since $m\geq 2$ and $\Phi_0(x)=x$ for every $x$, we have $f_0(x)=D^m\Phi_0(x)=0$. By uniqueness of the solution, it follows that $f_t(x)=D^m\Phi_t(x)$ for every positive $t$. We can now make estimates on $D^m\Phi_t(x)$ using this explicit expression. Applying \eqref{Fm} and \eqref{s estimate} gives
\begin{align*}\label{final lyapunov}
    \|D^m\Phi_t\|&\leq\int_0^t\|S(t,s,x)\|\cdot\|F_{m,s}\|ds\\
    &\leq c_m(1+\|X\|_m)\int_0^te^{-\mu(t-s)}e^{-2\mu s}ds\\
    &=c_m(1+\|X\|_m)e^{-\mu t}\int_0^te^{-\mu s}ds\\
    &\leq c_m(1+\|X\|_m)e^{-\mu t}.
\end{align*}
This concludes the proof.\qed
\bibliography{bibliography}

@article {Hamilton,
    AUTHOR = {Hamilton, Richard S.},
     TITLE = {The inverse function theorem of {N}ash and {M}oser},
   JOURNAL = {Bull. Amer. Math. Soc. (N.S.)},
  FJOURNAL = {American Mathematical Society. Bulletin. New Series},
    VOLUME = {7},
      YEAR = {1982},
    NUMBER = {1},
     PAGES = {65--222},
      ISSN = {0273-0979,1088-9485},
   MRCLASS = {58C15 (46G05 58D15)},
  MRNUMBER = {656198},
MRREVIEWER = {Peter\ Michor},
       DOI = {10.1090/S0273-0979-1982-15004-2},
       URL = {https://doi.org/10.1090/S0273-0979-1982-15004-2},
}

@article {Uhlenbeck,
    AUTHOR = {Uhlenbeck, K.},
     TITLE = {Generic properties of eigenfunctions},
   JOURNAL = {Amer. J. Math.},
  FJOURNAL = {American Journal of Mathematics},
    VOLUME = {98},
      YEAR = {1976},
    NUMBER = {4},
     PAGES = {1059--1078},
      ISSN = {0002-9327,1080-6377},
   MRCLASS = {58G99 (58C25)},
  MRNUMBER = {464332},
MRREVIEWER = {A.\ J.\ Tromba},
       DOI = {10.2307/2374041},
       URL = {https://doi.org/10.2307/2374041},
}

@book {Galdi,
    AUTHOR = {Galdi, G. P.},
     TITLE = {An introduction to the mathematical theory of the
              {N}avier-{S}tokes equations},
    SERIES = {Springer Monographs in Mathematics},
   EDITION = {Second},
      NOTE = {Steady-state problems},
 PUBLISHER = {Springer, New York},
      YEAR = {2011},
     PAGES = {xiv+1018},
      ISBN = {978-0-387-09619-3},
   MRCLASS = {35Q30 (35-02 76D03 76D05 76D07)},
  MRNUMBER = {2808162},
       DOI = {10.1007/978-0-387-09620-9},
       URL = {https://doi.org/10.1007/978-0-387-09620-9},
}

@book {transversalmapping,
    AUTHOR = {Abraham, Ralph and Robbin, Joel},
     TITLE = {Transversal mappings and flows},
      NOTE = {An appendix by Al Kelley},
 PUBLISHER = {W. A. Benjamin, Inc., New York-Amsterdam},
      YEAR = {1967},
     PAGES = {x+161},
   MRCLASS = {57.55},
  MRNUMBER = {240836},
MRREVIEWER = {J.\ Palis},
}

@book {Khalil,
    AUTHOR = {Khalil, Hassan K.},
     TITLE = {Nonlinear systems},
 PUBLISHER = {Macmillan Publishing Company, New York},
      YEAR = {1992},
     PAGES = {xii+564},
      ISBN = {0-02-363541-X},
   MRCLASS = {34-01 (34-02 34Cxx 34H05 93-01 93-02 93C10)},
  MRNUMBER = {1201326},
MRREVIEWER = {A.\ Halanay},
}

@article{Bando,
  title={GENERIC PROPERTIES OF THE EIGENVALUE OF THE LAPLACIAN FOR COMPACT RIEMANNIAN MANIFOLDS},
  author={Shigetoshi Bando and Hajime Urakawa},
  journal={Tohoku Mathematical Journal},
  year={1983},
  volume={35},
  pages={155-172},
  url={https://api.semanticscholar.org/CorpusID:122610348}
}

@book {MR1468735,
    AUTHOR = {Lee, John M.},
     TITLE = {Riemannian manifolds},
    SERIES = {Graduate Texts in Mathematics},
    VOLUME = {176},
      NOTE = {An introduction to curvature},
 PUBLISHER = {Springer-Verlag, New York},
      YEAR = {1997},
     PAGES = {xvi+224},
      ISBN = {0-387-98271-X},
   MRCLASS = {53-01 (53C21)},
  MRNUMBER = {1468735},
MRREVIEWER = {Man\ Chun\ Leung},
       DOI = {10.1007/b98852},
       URL = {https://doi.org/10.1007/b98852},
}

@incollection {NashOverview2,
    AUTHOR = {Andrews, Ben},
     TITLE = {Notes on the isometric embedding problem and the
              {N}ash-{M}oser implicit function theorem},
 BOOKTITLE = {Surveys in analysis and operator theory ({C}anberra, 2001)},
    SERIES = {Proc. Centre Math. Appl. Austral. Nat. Univ.},
    VOLUME = {40},
     PAGES = {157--208},
 PUBLISHER = {Austral. Nat. Univ., Canberra},
      YEAR = {2002},
      ISBN = {0-7315-5204-0},
   MRCLASS = {53C40 (58C15)},
  MRNUMBER = {1953483},
MRREVIEWER = {J\"urgen\ Eichhorn},
}

@article{Nash,
    AUTHOR = {Nash, John},
     TITLE = {The imbedding problem for {R}iemannian manifolds},
   JOURNAL = {Ann. of Math. (2)},
  FJOURNAL = {Annals of Mathematics. Second Series},
    VOLUME = {63},
      YEAR = {1956},
     PAGES = {20--63},
      ISSN = {0003-486X},
   MRCLASS = {53.1X},
  MRNUMBER = {75639},
MRREVIEWER = {J.\ Schwartz},
       DOI = {10.2307/1969989},
       URL = {https://doi.org/10.2307/1969989},
}

@article{NashOverview1,
  title={The masterpieces of John Forbes Nash Jr},
  author={De Lellis, Camillo},
  journal={arXiv preprint arXiv:1606.02551},
  year={2016}
}

@book {Srivastava,
    AUTHOR = {Srivastava, S. M.},
     TITLE = {A course on {B}orel sets},
    SERIES = {Graduate Texts in Mathematics},
    VOLUME = {180},
 PUBLISHER = {Springer-Verlag, New York},
      YEAR = {1998},
     PAGES = {xvi+261},
      ISBN = {0-387-98412-7},
   MRCLASS = {04A15 (28A05 54-01 54H05)},
  MRNUMBER = {1619545},
MRREVIEWER = {Marek\ Balcerzak},
       DOI = {10.1007/978-3-642-85473-6},
       URL = {https://doi.org/10.1007/978-3-642-85473-6},
}

@article{perelman2002entropy,
  title={The entropy formula for the Ricci flow and its geometric applications},
  author={Perelman, Grisha},
  journal={arXiv preprint math/0211159},
  year={2002}
}

@book {Evans,
    AUTHOR = {Evans, Lawrence C.},
     TITLE = {Partial differential equations},
    SERIES = {Graduate Studies in Mathematics},
    VOLUME = {19},
   EDITION = {Second},
 PUBLISHER = {American Mathematical Society, Providence, RI},
      YEAR = {2010},
     PAGES = {xxii+749},
      ISBN = {978-0-8218-4974-3},
   MRCLASS = {35-01},
  MRNUMBER = {2597943},
MRREVIEWER = {Diego\ M.\ Maldonado},
       DOI = {10.1090/gsm/019},
       URL = {https://doi.org/10.1090/gsm/019},
}

@article {Ambra,
    AUTHOR = {D'Ambra, Giuseppina and De Leo, Roberto and Loi, Andrea},
     TITLE = {Partially isometric immersions and free maps},
   JOURNAL = {Geom. Dedicata},
  FJOURNAL = {Geometriae Dedicata},
    VOLUME = {151},
      YEAR = {2011},
     PAGES = {79--95},
      ISSN = {0046-5755,1572-9168},
   MRCLASS = {53C42 (37J35 58A30)},
  MRNUMBER = {2780739},
MRREVIEWER = {Alma\ Luisa\ Albujer Brotons},
       DOI = {10.1007/s10711-010-9520-9},
       URL = {https://doi.org/10.1007/s10711-010-9520-9},
}

@book {Gilbarg,
    AUTHOR = {Gilbarg, David and Trudinger, Neil S.},
     TITLE = {Elliptic partial differential equations of second order},
    SERIES = {Classics in Mathematics},
      NOTE = {Reprint of the 1998 edition},
 PUBLISHER = {Springer-Verlag, Berlin},
      YEAR = {2001},
     PAGES = {xiv+517},
      ISBN = {3-540-41160-7},
   MRCLASS = {35-02 (35Jxx)},
  MRNUMBER = {1814364},
}

@book {MR681859,
    AUTHOR = {Aubin, Thierry},
     TITLE = {Nonlinear analysis on manifolds. {M}onge-{A}mp\`ere equations},
    SERIES = {Grundlehren der mathematischen Wissenschaften [Fundamental
              Principles of Mathematical Sciences]},
    VOLUME = {252},
 PUBLISHER = {Springer-Verlag, New York},
      YEAR = {1982},
     PAGES = {xii+204},
      ISBN = {0-387-90704-1},
   MRCLASS = {58-02 (35J60 35N99 53-02 58G30)},
  MRNUMBER = {681859},
MRREVIEWER = {A.\ G.\ Ramm},
       DOI = {10.1007/978-1-4612-5734-9},
       URL = {https://doi.org/10.1007/978-1-4612-5734-9},
}
\bibliographystyle{acm}
\end{document}